\documentclass[pdflatex,sn-mathphys-num]{sn-jnl}

\usepackage{amsmath,amssymb,amsfonts}
\usepackage{amsthm}
\usepackage{latexsym}
\usepackage{graphicx}
\usepackage{wrapfig}

\usepackage[english]{babel}
\addto\captionsenglish{\renewcommand{\figurename}{Fig.}}
\usepackage[utf8]{inputenc}
\usepackage[T1]{fontenc}
\usepackage{mathrsfs}
\usepackage{enumerate}
\usepackage{version}
\usepackage{calc}
\usepackage{subfigure}
\usepackage{tikz}
\usepackage{xcolor}
 \usepackage{float}
 \usepackage{indentfirst}
 \usepackage{hyperref}
 \usepackage{cleveref}
\hypersetup{hidelinks}
\usepackage{etoolbox}
\makeatletter
\patchcmd{\@maketitle}{\vskip21pt}{\vskip0pt}{}{}
\patchcmd{\@maketitle}{\vskip20pt}{\vskip10pt}{}{}
\patchcmd{\@maketitle}{\vskip20pt}{\vskip12pt}{}{}
\patchcmd{\@maketitle}{\vskip24pt}{\vskip10pt}{}{}
\patchcmd{\@maketitle}{\vskip24pt}{\vskip10pt}{}{}
\patchcmd{\@maketitle}{\vskip36pt}{\vskip12pt}{}{}
\patchcmd{\@maketitle}{Corresponding author(s). E-mail(s):}{Oleksandr Misiats (corresponding author):}{}{\PackageError{manuscript}{Could not update corresponding-author label}{}}
\patchcmd{\@maketitle}{Contributing authors:}{Sergio Conti:}{}{\PackageError{manuscript}{Could not update author email label}{}}
\renewcommand{\fnum@figure}{\figurename\ \thefigure}
\makeatother

\newcommand{\R}{{\mathbb R}}

\newcommand{\N}{{\mathbb N}}

\renewcommand{\d}{\partial}
\newcommand{\dist}{{\rm dist}\,}

\newcommand{\sym}{{\rm sym}\,}

\newcommand{\eps}{\varepsilon}
\newcommand{\ve}{\varepsilon}

\newcommand{\bp}{{\bf{p}}}
\newcommand{\bq}{{\bf{q}}}
\theoremstyle{thmstyleone}
\newtheorem{theorem}{Theorem}[section]
\numberwithin{equation}{section}
\newtheorem{lemma}[theorem]{Lemma}
\newtheorem{proposition}[theorem]{Proposition}

\newcommand\KM{\mathrm{KM}}
\newcommand\CC{\mathrm{CC}}

\newcommand{\LM}[1]{\hbox{\vrule width.2pt \vbox to#1pt{\vfill \hrule width#1pt height.2pt}}}
\newcommand{\LL}{{\mathchoice
{\,\LM7\,}{\,\LM7\,}{\,\LM5\,}{\,\LM{3.35}\,}}}

\newcommand\Es{E^\mathrm{scal}}
\newcommand\Eseps{E^\mathrm{scal}_\ve}

\newcommand\ECCone{E_1^{\CC}}
\newcommand\ECCeps{E_\eps^{\CC}}

\renewcommand{\nabla}{D}

\renewcommand\theenumi{{[\roman{enumi}]}}

\newcommand{\calH}{\mathcal H}

\begin{document}

\title{Scaling laws for linear-elastic two-variant martensitic phase transitions in general domains}

\author[1]{\fnm{Sergio} \sur{Conti}}\email{sergio.conti@uni-bonn.de}
\author*[2]{\fnm{Oleksandr} \sur{Misiats}}\email{omisiats@vcu.edu}
\affil[1]{\orgdiv{Institute for Applied Mathematics},
\orgname{University of Bonn}, \orgaddress{\city{Bonn}, \postcode{53115}, \country{Germany}}}
\affil[2]{\orgdiv{Department of Mathematics and Statistics},
\orgname{Virginia Commonwealth University}, \orgaddress{\city{Richmond}, \state{VA}, \postcode{23284}, \country{USA}}}

\abstract{
We consider a two-well problem within geometrically linear elasticity, on a bounded Lipschitz set in dimension $n\ge 2$. We assume that the two wells are
compatible through a symmetrized rank-one connection and that the displacement
matches their average on the boundary, and we include a regularization based on the total surface energy.
We show that there are two distinct regimes,
with the energy scaling differently in terms of the surface-energy coefficient $\eps$: it is proportional to $\eps^{2/3}$ if there are two distinct compatible twin-normal directions, and to $\eps^{4/5}$ if there is only one. This matches previous findings in rectangles and other special geometries, mostly in two dimensions, and the lower bounds obtained by a Fourier method by Rai\c{t}\u{a}, R\"uland, Tissot and Tribuzio (2024). Our results are obtained by explicit constructions, which contain a self-similar refinement of oscillations close to the boundary, and by a short real-space proof of the matching lower bounds.
}


\maketitle

\section{Introduction} 

A {\it martensitic phase transition} is a solid-to-solid phase transformation in certain metallic alloys, involving an abrupt change of the lattice  structure at a critical temperature. 
At high temperature the material is in a high-symmetry state; at low temperature, two or more variants of a lower-symmetry lattice structure appear. In most situations a fine mixture of different variants is observed, with very characteristic geometric structures.
The resulting variants of the deformation of the lattice unit cell are called {\it variants of martensite}. The transformation involves coordinated atomic displacements and a change of lattice geometry, without atomic diffusion.
The formation and rearrangement of these patterns are central to the shape-memory effect and to the mechanical response of transforming alloys; see Bhattacharya~\cite{Bha03}. A basic mathematical question is how the geometry of the transformation strains, the boundary conditions, and the cost of interfaces together determine the energy and length scales of the resulting microstructures.

The variational theory of martensite describes these phenomena through nonconvex elastic energies whose wells represent stress-free configurations. The foundational work of Ball and James~\cite{BJ87,BJ92} established a framework in which fine phase mixtures arise as energy-minimizing sequences, with compatibility restricting which mixtures can be realized by deformations. In geometrically linear elasticity, the corresponding relaxation problem was studied by Kohn~\cite{Koh91}. These developments connect crystallographic microstructures with quasiconvexity, Young measures, and differential constraints; a systematic account is given by M\"uller~\cite{Mul99}. Relaxation describes the macroscopic response of arbitrarily fine mixtures, but it does not by itself select a positive microscopic length scale. Such a scale emerges when elastic energy is supplemented by an interfacial or higher-gradient contribution.
As customary, we denote by
$\eps$ the energy per unit area associated with these interfaces.

The role of this regularization was demonstrated by Kohn and M\"uller~\cite{KM92,KM94} in their study of twinning near an austenite--martensite interface. In their scalar model, the elastic energy favors alternating regions in which a derivative takes one of two preferred values, while the boundary data require their oscillations to disappear at the interface. A laminate of uniform width can satisfy this condition only at an additional elastic cost
close to the interface, leading to
a minimum energy of order $\varepsilon^{1/2}$. A branching pattern instead refines the twins progressively toward the interface, balancing elastic accommodation against the increasing number of interfaces. The resulting minimum energy scales as $\varepsilon^{2/3}$, with matching upper and lower bounds. Local energy estimates and asymptotic self-similarity properties of minimizers have  subsequently been established in~\cite{Con00}. More recently, refined local bounds and self-similarity results for various boundary
conditions, arbitrary rectangular domains, and arbitrary volume fractions of the martensitic
variants have been obtained~\cite{CDKZ21}.

Branching is not specific to martensitic transformations. Related competitions between fine-scale accommodation and interfacial cost 
leading to branching patterns
occur, for example, in uniaxial ferromagnets~\cite{CKO99}, 
 diblock copolymers~\cite{Choksi01},
the Ginzburg--Landau model of type-I superconductors~\cite{COS16}, 
folding patterns in thin films~\cite{BCDM02,bella-kohn:14,BourneContiMueller2017}, 
optimal design~\cite{KW14,KW16}
and other systems. These examples illustrate the broader usefulness of matching upper bounds obtained from explicit constructions with ansatz-free lower bounds. They also underline that the differential constraint and the geometry of the preferred states, rather than branching alone, determine the relevant scaling law; we refer to \cite{Koh07} for an overview.

Passing from a scalar model to elasticity introduces the additional difficulty that the elastic energy depends only on the symmetric part of the displacement gradient and therefore does not directly penalize infinitesimal rotations, so that compatibility and rigidity must be considered together with interfacial regularity. Structural aspects of such vector-valued models have been studied for two-well microstructures with finite surface energy~\cite{DM95}, in a sharp-interface limit for geometrically linear elasticity~\cite{CS06},  through quantitative rigidity estimates for the cubic-to-tetragonal transformation~\cite{CO12},
and for the rigidity of branching microstructures~\cite{Sim21}.
Branching models have also been extended to multiparameter settings: scaling regimes and optimal needle shapes motivated by low-hysteresis shape-memory alloys~\cite{Zwi14}, partially branched regimes at low volume fraction~\cite{CZ16}, as well as geometry- and parameter-dependent matching bounds for a martensitic nucleus embedded in an elastic austenitic matrix~\cite{CDMZ20}. These results emphasize that, in vector-valued and multiparameter problems, branching scaling laws must be understood together with rigidity constraints, and a fixed-parameter power law should be distinguished from a uniform description of the full parameter-dependent phase diagram.

A different mechanism arises when the geometry of the two wells becomes degenerate. Chan and Conti~\cite{CC14,CC15} studied two-dimensional nonlinear elastic models with $SO(2)$ invariance and showed that the number of rank-one connections can change the branching exponent.
 Their analysis distinguishes a regime with two rank-one connections and $\varepsilon^{2/3}$ scaling from a single-connection regime with $\varepsilon^{4/5}$ scaling. The latter exploits a vectorial cancellation that is absent from a scalar reduction: one displacement component compensates a shear strain generated by the other. 
Although one might expect that fewer compatible interface directions lead to a larger minimum energy, the interplay with invariance under rotations leads to a smaller one. Let us address this in the geometrically linear framework used in this paper (the nonlinear setting is almost the same). Rank-one connections in two dimensions 
between two wells determined by matrices $A$ and $B$ can be  seen as zeros of the determinant 
of $A-B-S$, with $S$ skew symmetric; the determinant is a quadratic expression. If $A$ and $B$ are such that there is only one solution, then it is a degenerate solution, and $|\det (A-B-S)|$ increases only quadratically as one moves away from it. In other words, for neighbouring values of $S$
the matrix $A-B-S$ is also approximately rank-one, so that slightly tilted interfaces have very low elastic energy. Therefore it is possible to bend the domains and obtain a more efficient branching construction. In the case of two rank-one connections, instead, the quadratic equation has two simple solutions.
One is used in the construction, and having a second one at a far-away point does not influence the scaling of the energy.

Twinning with spatially varying volume fraction provides another setting in which this mechanism matters.
 The work~\cite{CKM22} established scaling laws for two nonlinear elastic phases with a single rank-one connection under bending-type boundary conditions, obtaining the scaling $\eps^{3/4}$ or $\eps^{4/5}$,
depending on details of the boundary conditions.  In~\cite{CKM24} bending-induced twinning for a three-dimensional geometrically linear model with two martensitic variants was studied using displacement or loading conditions that enforce bending, leading to the  scaling $\eps^{2/3}$.

Two recent papers are particularly relevant for this work, even if we do not use their results. They place such scaling laws in the more general setting of fields satisfying linear differential constraints. Rai\c{t}\u{a}, R\"uland, and Tissot~\cite{RRT23}
established $\varepsilon^{2/3}$ lower bounds for classes of $\mathcal A$-free fields
on Lipschitz domains for first-order operators.
Rai\c{t}\u{a}, R\"uland, Tissot, and Tribuzio~\cite{RRTT24} extended this approach to higher-order 
operators, with exponent $2L/(2L+1)$, where $L$  is a measure of degeneracy of the compatible interfaces, analogous to the order of the zero of $\det(A-B-S)$ in the discussion above. In practice, $L=1$ and $L=2$ in the two cases of interest here, so that
for the symmetrized gradient they obtain lower bounds
of order $\varepsilon^{2/3}$ and $\varepsilon^{4/5}$ on bounded Lipschitz domains in any dimension. Matching upper bounds are provided for a square in two dimensions.
They use a different regularization, which is equivalent at the level of scaling, see~R\"uland, Tissot, Tribuzio and Zillinger (2025, preprint, \href{https://arxiv.org/abs/2507.06773}{arXiv:2507.06773}).
A general-domain result therefore primarily requires constructing a pattern that realizes the predicted scaling under the full boundary constraint.

 Staircase-type configurations in linearized elasticity
 have been investigated in~\cite{MR26},
showing that lower scaling bounds depend both on
lamination order and on the occurrence of degenerate symmetrized rank-one directions. Their constructions in selected geometries extend the role of this degeneracy beyond first-order laminates.
Incompatible two-well problems, including linearized elasticity, have been studied in~a preprint by Piemontese-Fischer (2025, \href{https://arxiv.org/abs/2512.12014}{arXiv:2512.12014}) after separating off the relaxed bulk-energy contribution. In the two-dimensional elasticity results, general-domain lower bounds are complemented by matching constructions on suitably oriented squares.
Domain dependence was recently investigated in~\cite{GRTZ26} for a related class of problems describing martensitic inclusions in an austenitic matrix.
In  special domains exact solutions with finite boundary energy exist, and they clearly have energy of order $\eps$. Their main result is that for generic domains the energy instead has the higher
 $\eps|\ln\eps|$ scaling. 
 This is related to nucleation problems, in which the shape 
of the
martensitic region is itself an unknown~\cite{KK11,KKO13}, and in which corners and boundaries
can lower the energy barrier~\cite{BG15,TZ25}. 

With the exception of the domain-geometry results in~\cite{GRTZ26}, 
and of the lower bounds in~\cite{RRT23,RT24,RRTT24},
essentially all of
the sharp vectorial and geometrically linear scaling laws mentioned in the previous paragraphs have so far been established only in special
two-dimensional domains: rectangles or parallelograms with sides aligned with the
relevant rank-one directions, and occasionally cubes in higher dimension~\cite{CZ16,RT24,MR26}. The purpose
of the present paper is to remove this restriction,
both for the
$2/3$ and for the $4/5$ exponents.

\subsection{Setting of the problem and main result}
Let $\Omega\subset\mathbb R^n$, $n\geq2$, be a nonempty bounded Lipschitz open set, and let $\bp$, $\bq\in\mathbb R^n$ be nonzero. We write
\[
 e(u):=\frac{\nabla u+\nabla u^T}{2},
 \qquad
 (\bp\otimes \bq)_{\mathrm{sym}}
 :=\frac{\bp\otimes \bq+\bq\otimes \bp}{2}.
\]
For definiteness we work with homogeneous boundary conditions, which correspond to the average of the two eigenstrains $\pm (\bp\otimes \bq)_\sym$. For
\begin{equation}\label{eqfunctionspace}
 u\in W^{1,2}_0(\Omega;\mathbb R^n),
 \qquad \nabla u\in BV(\Omega;\mathbb R^{n\times n}),
\end{equation}
we consider
\begin{equation}\label{eq:energy}
\begin{split}
 E_\varepsilon[u,\Omega]
 := {}&\int_\Omega\min\left\{
 \left|e(u)-(\bp\otimes \bq)_{\mathrm{sym}}\right|^2,
 \left|e(u)+(\bp\otimes \bq)_{\mathrm{sym}}\right|^2
 \right\}\,dx
 +\varepsilon |D^2u|(\Omega).
\end{split}
\end{equation}
The first term measures departure from the two preferred strains, while the second penalizes the total variation of the displacement gradient. The boundary condition is an idealization of a surrounding medium that prescribes the macroscopic displacement; it does not include the elastic response of an austenitic matrix as an additional variable.

The representation of the wells makes their compatibility explicit. If $\bp$ and $\bq$ are not parallel, the two compatible twin-normal directions are those of $\bp$ and $\bq$, up to normalization and reversal. If they are parallel, these directions coincide and the difference of the preferred strains has rank one. Our main result shows that this distinction determines the exponent for every fixed domain in the above class.

\begin{theorem}[Energy scaling]\label{thm:main}
Let $\Omega\subset\mathbb R^n$, $n\ge2$, be a nonempty bounded Lipschitz open set, and let $\bp$, $\bq\in\mathbb R^n$ be nonzero. There exists $C=C(\Omega,\bp,\bq)>0$ such that
for any $\varepsilon\in(0,1]$,
\[
\begin{array}{ll}
 C^{-1}\varepsilon^{2/3}
 \leq \displaystyle\inf_u E_\varepsilon[u,\Omega]
 \leq C\varepsilon^{2/3},
 &\text{if $\bp$ and $\bq$ are not parallel},\\[7pt]
 C^{-1}\varepsilon^{4/5}
 \leq \displaystyle\inf_u E_\varepsilon[u,\Omega]
 \leq C\varepsilon^{4/5},
 &\text{if $\bp$ and $\bq$ are parallel}.
\end{array}
\]
The infima are taken over the admissible displacements specified in~\eqref{eqfunctionspace}.
\end{theorem}

Let us briefly explain our strategy, first for the  upper and then for the lower bound. We intend to construct a map $u\in W^{1,2}_0(\Omega;\R^n)$ which oscillates in the direction $\bq$, with domain branching close to the boundary. This involves picking $\bp\otimes \bq$ as the relevant rank-one connection. The branching direction cannot be fixed a priori: the typical thickness of the laminate around a point $x$ has to go to zero as $x$ approaches $\partial\Omega$. Therefore the refinement has to follow the boundary; we shall use the distance to the boundary to make this quantitative.
In order to obtain $u$, we start from a construction 
$v:[0,\infty)\times \R \to\R$ for the scalar version of the problem defined on $\R^2$, with 
$\partial_2 v$ close to $\pm 1$, and $v(0,\cdot)=0$.
The first attempt is to 
 use 
 \begin{equation}\label{u:rescaled}
 u(x):=\bp v(\dist(x,\partial\Omega), \bq\cdot x).
  \end{equation}
  Then $Du=\bp\otimes \bq\partial_2 v + \bp\otimes \nabla\dist(\cdot,\partial\Omega)\partial_1 v$ is close to $\pm \bp\otimes \bq$, with additional errors coming from $\partial_1 v$, and $u$ vanishes on $\partial\Omega$.
Integrals over $\Omega$ are then turned into integrals over a two-dimensional domain using the coarea formula, if one can bound the Jacobian of the relevant map, which depends on both the modulus of  $\nabla \dist(\cdot,\partial\Omega)$ and its orientation with respect to $\bq$.
 This strategy has several difficulties. First, computing second derivatives of $u$ is problematic, since the function $x\mapsto \dist(x,\partial\Omega)$ does not, in general, have classical second derivatives. In smooth domains, it is smooth close to the boundary, but not in the interior. In convex polyhedra its gradient $\nabla \dist(\cdot,\partial\Omega)$ is piecewise constant, and jumps on hypersurfaces of codimension 1. In particular, it is a function of bounded variation. Therefore~\eqref{u:rescaled} can only be used directly in a class of convex  polyhedra; in a general Lipschitz domain the distance does not possess classical second derivatives, and its distributional derivatives might interact in a complex way with the distributional derivatives of $v$.
 
 A second problem arises in the case that $\bp$ and $\bq$ are parallel. Here, even for simple domains and even staying away from singularities of $\nabla\dist(\cdot,\partial\Omega)$, the expression in~\eqref{u:rescaled}
 does not have the optimal scaling. After reducing  without loss of generality to the case  $\bp=\bq$, we show below that in regions where $\dist(\cdot,\partial\Omega)$ is smooth we can use,
starting from the two-dimensional construction
$\tilde v:[0,\infty)\times \R \to\R^2$
relative to the wells $\pm e_2\otimes e_2$,
 \begin{equation}\label{u:rescaled45}
\tilde u(x):=\bp \tilde  v_2(\dist(x,\partial\Omega), \bq\cdot x)
+\nabla \dist(x,\partial\Omega)  \tilde v_1(\dist(x,\partial\Omega), \bq\cdot x).
  \end{equation}
 This expression exploits cancellation between the first-order derivatives $\partial_1\tilde v_2$ and $\partial_2\tilde v_1$ in the shear strain. This cancellation is the origin of the difference between the two scaling regimes in the theorem above, and therefore needs to be preserved when changing the geometry.
 However, the expression in \eqref{u:rescaled45} is even less regular than the one in \eqref{u:rescaled}, as it already contains derivatives of the distance function, which are  typically discontinuous.  Therefore one needs a regularized version of the distance, which still vanishes on the boundary, and still allows one to use the coarea formula. 
In turn, in domains with corners the regularized distance function necessarily has second
derivatives which diverge close to the boundary, requiring a more careful treatment. In order to make the construction independent of details of the domain, for both scalings we first obtain an explicit construction in a cube, in the spirit of~\cite{CZ16,RT23}, and then show that a general domain can be covered by countably many cubes, with an appropriate bound on the energy (in practice, we shall need a boundary layer of thickness $\eps$ so that the number of effectively used cubes is finite and of order $\eps^{1-n}$).

It is worth noting that the construction used for the upper bound in general domains introduces a number of additional internal interfaces, which in the end do not affect the scaling of the energy.
At the same time, the strategy used for cubes generalizes to a class of simple domains; for example, in the $2/3$ case one can easily prove it for generic convex polyhedra. It is easy to imagine that many more domains can be treated that way, generating patterns that are closer to what one expects in reality,  but for brevity, we do not discuss the details.

Matching lower bounds have been obtained by Fourier methods in~\cite{RRTT24}. We give here a different proof in real space 
using the localization method in the
spirit of \cite{CC15}: small cubes are classified according to whether
the surface energy on that cube is large, the strain deviates substantially from the
wells, or the displacement itself is large, without relying on any particular ansatz
for the minimizer. When $\bp$ and $\bq$ are not parallel, directions orthogonal to at least one of these vectors span $\mathbb R^n$.
Along such directions, the two-well structure of the energy disappears, and directional Poincar\'e estimates control all displacement components by the elastic energy. In the parallel case, the same argument leaves the displacement component along the common direction uncontrolled. {Recovering control uses a null-Lagrangian identity for $2\times2$ minors, followed by a transverse Poincar\'e inequality.} This loss of direct coercivity is the lower-bound counterpart of the cancellation exploited by the upper-bound construction.

The paper is organized as follows. Section~\ref{secupperbound}  develops the upper bounds, treating the nonparallel and parallel cases separately and extending the planar constructions first to cubes and then to general domains. Section~\ref{seclowerbound} is devoted to the lower bounds.
We first present some general localization tools, then explain the localization method for a full-gradient two-well energy. The proof is then presented separately for the symmetric-gradient problem and the case with an additional degeneracy when $\bp$ and $\bq$ are parallel.

\section{Construction for the upper bound }\label{secupperbound}

\subsection{\texorpdfstring{Upper bound when $\bp$ and $\bq$ are not parallel}{Upper bound when p and q are not parallel}}
\subsubsection{The Kohn--M\"uller construction for a half-plane}
\label{secKMhalfp}
We start by reviewing the by-now classical branching construction by Kohn and Müller \cite{KM92,KM94} (Theorem~\ref{Th:KM} below). For what follows, we need some additional bounds on the function and its derivatives, which are obtained in Lemma~\ref{lem:KMbounds} below.

For $\varepsilon>0$, $\omega\subset\R^2$ open and bounded, and $u\in W^{1,2}(\omega)$ with $\nabla u\in BV(\omega;\R^2)$, we write
\begin{equation*}
E^\varepsilon_\KM[u, \omega]
 :=  \int_{\omega} ( |\d_1 u|^2 + \min\{(\d_2 u \pm 1)^2\} + \ve |D^2 u|) \, dx dy
\end{equation*}
where, as usual in this context, the integral of the last term represents the total variation of the measure $D^2u$.
In the entire paper we write $D=(\d_1, \dots, \d_n)$ for both the distributional and the weak gradient.
In this section we use $(x,y)$ for points in $\R^2$.

\newcommand\widetildeuKM{u_\KM}
\begin{theorem}\label{Th:KM}
There is $C>0$ such that for any
 $\ve>0$ there is a function $\widetildeuKM\in W^{1,2}_\mathrm{loc}((0,\infty) \times \R)$
which satisfies the following properties: 
\begin{enumerate}
    \item\label{Th:KMzero} $\widetildeuKM(0,y) = 0$ for all $y\in\R$, in the sense of traces;
    \item\label{Th:KMx23} $|\widetildeuKM(x,y)| \leq C \ve^{1/3} x^{2/3}$ pointwise;
    \item\label{Th:KMy1} $|\partial_2\widetildeuKM(x,y)| \le 1$ almost everywhere;
    \item\label{Th:KMenergy}
    $D\widetildeuKM\in BV_\mathrm{loc}((0,\infty)\times \R;\R^2)$ and for any $L>\eps$,
    $H>0$, $a \in \R$ we have
    \begin{equation}\label{KM energy}
E^\varepsilon_\KM[\widetildeuKM, (\eps,L)\times (a,a+H)]
\leq C (\ve^{2/3}L^{1/3}H  + \ve L).
 \end{equation}
 \end{enumerate}   
\end{theorem}
We remark that the model in~\cite{KM94} contained only $\eps|\partial_2\partial_2 u|$
as regularization, and not {the full Hessian $D^2u$. This allowed the authors to obtain an} energy estimate like~\eqref{KM energy} in the full domain $(0,L)\times (a,a+H)$, while keeping $\partial_2 u=\pm1$ pointwise almost everywhere. Due to the change of variables presented in~\eqref{u:rescaled}, we need to estimate the full Hessian $D^2u$ (see~\eqref{second_der} below), and therefore use a different treatment  of the boundary layer, which will be introduced in
Lemma~\ref{lem:KMbounds} below.
\begin{proof}
Assume for a moment that we can prove the theorem for $\eps=1$, obtaining a constant $\hat C$ and a function $\hat u$ as stated. We then define $C:=\hat C$ and
\begin{equation}\label{eqresacaleu}
\widetildeuKM(x, y):=\eps \hat u(x/\eps , y/\eps).
\end{equation}
Properties~\ref{Th:KMzero}, \ref{Th:KMx23}  and~\ref{Th:KMy1} are immediately inherited.
To obtain~\ref{Th:KMenergy} we compute
$D\widetildeuKM=D\hat u$, $D^2\widetildeuKM=\eps^{-1}D^2\hat u$ (evaluated at the appropriate points)  and then, with a change of variables,
\begin{equation}
\begin{aligned}
 &E^\varepsilon_\KM[\widetildeuKM,(\eps,L)\times(a,a+H)]\\
 &\quad=\eps^2\int_1^{L/\eps}\int_{a/\eps}^{(a+H)/\eps}
 \bigl(|\d_1\hat u|^2+\min\{(\d_2\hat u\pm1)^2\}+|D^2\hat u|\bigr)\,dy\,dx\\
 &\quad
 \le\eps^2\hat C\bigl((L/\eps)^{1/3}(H/\eps)+L/\eps\bigr)
 =\hat C\bigl(\eps^{2/3}L^{1/3}H+\eps L\bigr).
\end{aligned}
\end{equation}

For the remainder of the proof, we take $\eps=1$.
The details of the proof are provided in \cite{KM94}. For completeness and later reference, we provide a sketch of the construction. Fix
\begin{equation}\label{eqchoicetheta23}
 \theta:=2^{-3/2}\in (\frac14,\frac12), \hskip5mm \text{ so that }
\theta^{1/3}=2\theta
\end{equation}
and, for any $i\in\mathbb Z$, set
 $x_i :=  \theta^i$. The construction  consists of three steps:\\
{\bf Step 1.} At the heart of the construction in \cite{KM94} is the period-doubling cell problem, described by the function $v:[0,1]\times \R\to\mathbb{R}$
defined for $y\in [0,1/2)$ by
\begin{equation}\label{eqdefv}
v(x,y):=
\begin{cases}
y,
& \text{if } 0\leq y\leq \dfrac18+\dfrac{x}{8},
\\[6pt]
\dfrac14+\dfrac{x}{4}-y,
& \text{if } \dfrac18+\dfrac{x}{8}
   \leq y\leq \dfrac38+\dfrac{x}{8},
\\[6pt]
-\dfrac12+y,
& \text{if } \dfrac38+\dfrac{x}{8}
   \leq y< \dfrac12,
\end{cases}
\end{equation}
and then extended by
$v(x,y)=v(x,y+1)=-v(x,-y)$, see Figure~\ref{figvkm}.
It obeys $v(0,y)=\frac12v(1,2y)$ and $\partial_2v\in\{\pm 1\}$ almost everywhere.

\begin{figure}[t]
\centering
\begin{tikzpicture}[scale=.9]\small
\draw (0,0) rectangle (6,7);

\draw (0,0.875) -- (6,1.75);

\draw (0,2.625) -- (6,3.5);

\draw (0,4.375) -- (6,3.5);

\draw (0,6.125) -- (6,5.25);

\node at (3.1,0.65) {$\partial_2 v=1$};
\node at (3.1,2.15) {$\partial_2 v=-1$};
\node at (3.1,3.65) {$\partial_2 v=1$};
\node at (3.1,4.95) {$\partial_2 v=-1$};
\node at (3.1,6.45) {$\partial_2 v=1$};

\node[below] at (3,0) {$v(x,0)=0$};
\node[above] at (3,7) {$v(x,1)=0$};

\node[left] at (0,0.875) {$v(0,1/8)=1/8$};
\node[left] at (0,2.625) {$v(0,3/8)=-1/8$};
\node[left] at (0,4.375) {$v(0,5/8)=1/8$};
\node[left] at (0,6.125) {$v(0,7/8)=-1/8$};

\node[right] at (6,1.75) {$v(1,1/4)=1/4$};
\node[right] at (6,5.25) {$v(1,3/4)=-1/4$};

\foreach \y in {0, 0.25,.5, 0.75, 1} {
\draw (0,7*\y) circle (0.09); }
\foreach \y in {0, .5, 1} {
\draw (6,7*\y) circle (0.09); }

\end{tikzpicture}
\caption{Sketch of the function $v$ used in the Kohn--M\"uller construction in the domain $[0,1]\times [0,1]$, see~\eqref{eqdefv}.
The circles mark the points at $x=0$ and $x=1$ where $v$ vanishes}
\label{figvkm}
\end{figure}
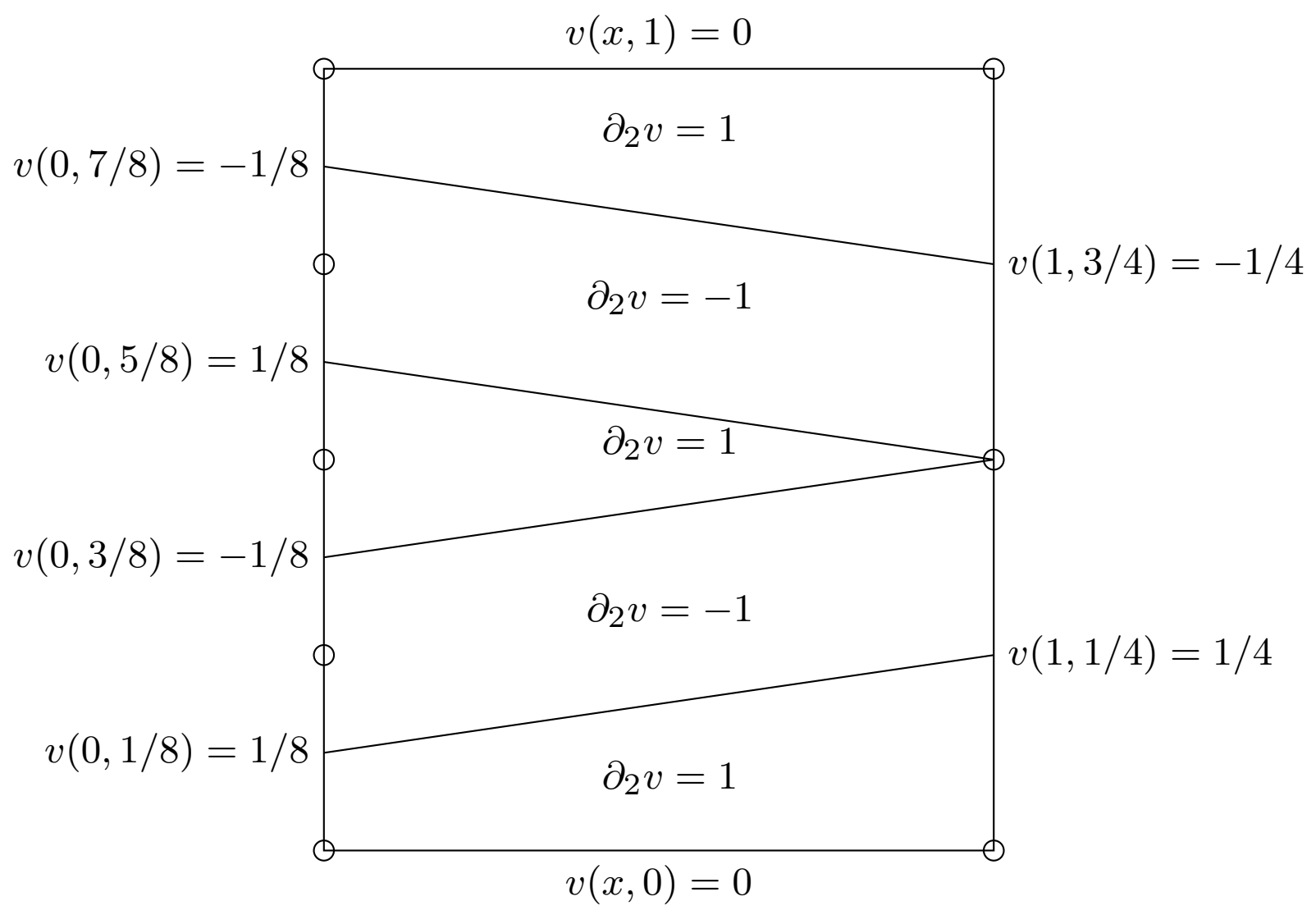 

{\bf Step 2.} For $x \in [x_1, x_0]$, the authors set
\[
\widetildeuKM(x,y)
:=
v\left(
\frac{x-x_1}{x_0-x_1},
y
\right).
\]
The resulting $\widetildeuKM$ satisfies
\[
|\partial_2\widetildeuKM|=1 \qquad \text{a.e.},
\]
\[
\widetildeuKM(x,y+1)=\widetildeuKM(x,y)
\qquad \text{for all } x,y,
\]
\[
\widetildeuKM(x_1,y)=\frac{1}{2} \widetildeuKM(x_0,2y)
\qquad \text{for all } y.
\]
{\bf Step 3.} The function $\widetildeuKM$ is extended to $(0,\infty)\times\mathbb{R}$ by an anisotropic rescaling,
\begin{equation}\label{eqdefutildeKM}
\widetildeuKM(x,y):=2^{-i} \widetildeuKM\left(\theta^{-i}x,2^i y\right)
=2^{-i}
v\left(
\frac{x-x_{i+1}}{x_i-x_{i+1}},
2^i y
\right)
\qquad \text{if } x\in[x_{i+1},x_i].
\end{equation}
Properties~\ref{Th:KMzero} and
\ref{Th:KMy1} are immediate.
Property~\ref{Th:KMx23} follows from $|\widetildeuKM(x,y)|\le 2^{-i} \|v\|_{L^\infty}$ for $x\in [x_{i+1},x_i]$ and
$\theta^{2/3}=1/2$ (see~\eqref{eqchoicetheta23}).
The energy estimate~\ref{Th:KMenergy} was established in \cite{KM94} for the first-order term and only for the $\d^2_2\widetildeuKM$ term in the second-order one;
the remaining two terms can be treated similarly, see for example~\cite{CDMZ20}.
For the convenience of the reader we repeat the computation for the full second-order term.
Fix $i\le 0$ and let
  $B_i:=(x_{i+1},x_i]\times(0,2^{-i}]$, one period of the $i$-th strip.
  Counting the four interfaces of $v$ in~\eqref{eqdefutildeKM}, and for $\d_1^2  u_\KM$ also the
  jump of $\d_1 u_\KM$ across the segment $\{x_i\}\times(0,2^{-i}]$, one obtains
  \[
  |\d_2^2 u_\KM|(B_i)\le C\theta^i,\qquad
  |\d_1\d_2 u_\KM|(B_i)\le C2^{-i},\qquad
  |\d_1^2 u_\KM|(B_i)\le C(4\theta)^{-i}.
  \]
  Since $\theta\in(\frac14,\frac12)$ and
  the bound in~\ref{Th:KMenergy} only includes $x\ge 1$, that is,
  $i\le0$, the first bound is the largest,
so $|D^2 u_\KM|(B_i)\le C\theta^i$, and summing over the periods and strips as for the $\d_2^2\widetildeuKM$ term concludes the proof.
\end{proof}
We next modify the construction from~\cite{KM94} to insert a boundary layer, so that the energy estimate can extend to $x=0$.
 \begin{lemma}\label{lem:KMbounds}
The constant $C$ and the  function
$u = u_\KM:(0,\infty) \times \R \to \R$
of
 Theorem \ref{Th:KM} can be chosen to additionally satisfy
\begin{enumerate}\setcounter{enumi}{4}
\item\label{lem:KMboundsdxu}
for any $x>0$ and $y \in \R$ we have
 \begin{equation}\label{KMGrad}
|\d_1 u_\KM(x,y)| \leq C \frac{\ve^{1/3}}{x^{1/3}} ;
\end{equation}
\item\label{lem:KMboundszero}  $u_\KM(x,y) = 0$ for $x  \in (0, \ve)$;
\item\label{enumcinfty}  $u_\KM\in C^\infty((0,\infty)\times \R)$;
\item\label{lem:KMboundsenergy}
for every $L$, $H>0$ and $a\in\R$,
\begin{equation}\label{eq:KM-full-rectangle}
 E^\eps_\KM[u_\KM,(0,L)\times(a,a+H)]
 \le C\bigl(\eps^{2/3}L^{1/3}H+\eps L\bigr).
\end{equation}
\end{enumerate}
\end{lemma}
\begin{proof}
As in Theorem~\ref{Th:KM}, the rescaling in~\eqref{eqresacaleu} shows that it suffices to consider $\eps=1$.
We denote by
$\widetilde{u}_\KM$ the function constructed in Theorem~\ref{Th:KM} for $\eps=1$.

 {\bf Step 1}: We show that
\begin{equation}\label{KM_grad_aux}
    |\d_1 \widetilde{u}_\KM(x,y)| \leq \frac{1}{x^{1/3}} \hskip5mm\text{ almost everywhere}.
\end{equation}
Indeed,
from~\eqref{eqdefv} it is clear that
$|\partial_1 v|\le1/4$ almost everywhere.
For $x\in(\theta^{i+1},\theta^i)$, the rescaling in~\eqref{eqdefutildeKM} gives
\[
 |\partial_1\widetilde u_\KM(x,y)|
 \le\frac{2^{-i}}{4(\theta^i-\theta^{i+1})}
 =\frac1{4(1-\theta)(2\theta)^i}\le x^{-1/3},
\]
because $(2\theta)^i=\theta^{i/3}\ge x^{1/3}$
and $4(1-\theta)\ge1$.
This proves~\eqref{KM_grad_aux}.

{\bf Step 2.} We will now modify $\widetilde u_\KM$ so that it additionally satisfies~\ref{lem:KMboundszero}. To this end, let $\eta \in C^{\infty}(\R;[0,1])$ be such that
(see Figure~\ref{fig:eta})
\begin{equation}\label{eqchoiceeta}
\text{
$\eta(x)=0$ for $x \leq 1+\frac18$, \quad
$\eta(x)=1$ for $x \geq 2-\frac18$, \quad
$0 \leq \eta'(x) \leq 2$, \quad
and $ |\eta''(x)|\leq 20$.}
\end{equation}
\begin{figure}[t]
\centering
\begin{tikzpicture}[scale=1.2]

\draw[->] (-0.2,0) -- (3.4,0) node[above] {$x$};
\draw[->] (0,-0.5) -- (0,1.6) node[right] {$\eta$};

\fill (1,0) circle (0.025);
\fill (2,0) circle (0.025);
\node[below] at (1,0) {$1$};
\node[below] at (2,0) {$2$};

\fill (0,1) circle (0.025);
\node[left] at (0,1) {$1$};

\draw[thick] (-0.2,0) -- (1.1,0);

\draw[thick]
(1.1,0)
.. controls (1.25,0) and (1.35,0.15) ..
(1.5,0.5)
.. controls (1.65,0.85) and (1.8,1) ..
(1.9,1);

\draw[thick] (1.9,1) -- (2.9,1);


\end{tikzpicture}
\caption{Sketch of the function $\eta$ defined in~\eqref{eqchoiceeta}}
\label{fig:eta}
\end{figure}
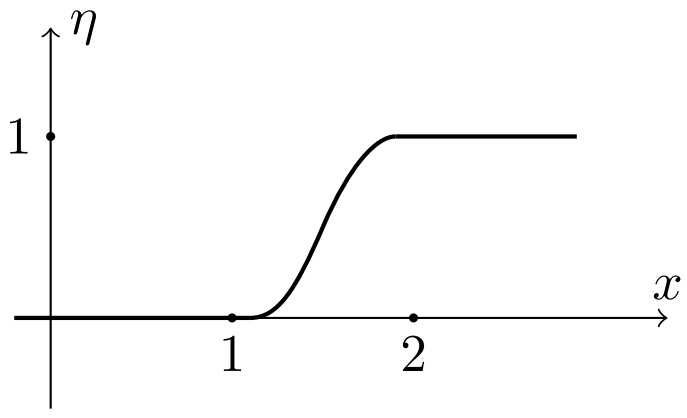
\newcommand\uKM{w}
We set
\[
\uKM(x,y):= \widetilde{u}_\KM(x,y) \eta\left(x\right).
\]
Since $|\eta| \leq 1$, it is straightforward to see that $\uKM$ satisfies~\ref{Th:KMzero}, \ref{Th:KMx23} and~\ref{Th:KMy1}
of Theorem \ref{Th:KM}.
Furthermore, since
\[
\d_1 \uKM(x,y) = \d_1 \widetilde{u}_\KM(x,y) \eta(x) + \widetilde{u}_\KM(x,y)  \eta'(x),
\]
the bound \ref{lem:KMboundsdxu} (for $\ve=1$) follows from  combining  Theorem \ref{Th:KM}\ref{Th:KMx23}
and \eqref{KM_grad_aux}.

{\bf Step 3.} It remains to verify the energy bound \eqref{eq:KM-full-rectangle}. Since {$\uKM=\widetilde u_\KM$} for $x\ge 2$,
and
$\uKM=0$ in the strip $x \in [0,1]$ while $\widetilde u_\KM - \uKM =\widetilde u_\KM(1 - \eta(x))$ in the strip $x \in [1, 2]$, we consider the two strips separately.

The estimate of the first strip is straightforward:
\begin{equation}\label{err1}\begin{split}
E^1_\KM[\uKM, (0,1)\times (a,a+H)]
 =& \int_{0}^{1}  \int_{a}^{a+H} 1 \, dy\, dx  =  H.
 \end{split}
\end{equation}
To proceed with the transition strip, since $0 \leq \eta \leq 1$, by~\ref{Th:KMy1} we have
\[
\int_{1}^{2}  \int_{a}^{a+H}  \min\{(\d_2 \uKM \pm 1)^2\} \,  dy\, dx\leq \int_{1}^{2}  \int_{a}^{a+H}  1 \,  dy \, dx =  H.
\]
Next, 
\begin{equation}\label{err2_1}
\int_{1}^{2 }  \int_{a}^{a+H} |\d_1 \uKM|^2 \, dy\, dx
\leq \int_{1}^{2}  \int_{a}^{a+H} [2 |\d_1 \widetilde u_\KM(x,y)|^2 + 2  \eta'(x)^2 \widetilde u_\KM^2(x,y)]\, dy\, dx\le CH
\end{equation}
where in the last step we used
property~\ref{Th:KMx23} from Theorem \ref{Th:KM} and  \eqref{KM_grad_aux}.
Finally, to estimate the surface energy contribution, we compute
\[
|D^2 \uKM| \leq  |\eta D^2 \widetilde u_\KM |  + 2 \eta'(|\d_2 \widetilde u_\KM|
+ |\d_1 \widetilde u_\KM|) + |\widetilde u_\KM \eta''|,
\]
which, in view of~\ref{Th:KMx23} of Theorem \ref{Th:KM} and \eqref{KM_grad_aux} yields
\begin{equation}\label{err2_3}
     \int_{1}^{2}  \int_{a}^{a+H} |D^2 \uKM| \, {dy\,dx} \leq   \int_{1}^{2}  \int_{a}^{a+H} |D^2 \widetilde u_\KM| \, {dy\,dx} + C  H.
\end{equation}

Combining~\eqref{err2_1}, \eqref{err2_3}, and
~\eqref{KM energy} for $\widetilde u_\KM$ gives
$E^1_\KM[\uKM,(1,2)\times(a,a+H)]\le C(H+1)$.
Together with~\eqref{err1} and the unchanged construction for $x\ge2$,
this proves~\eqref{eq:KM-full-rectangle} for $L\ge1$ and $\eps=1$.
For $0<L<1$, as in~\eqref{err1} one obtains the bound $LH\le L^{1/3}H$. This concludes the proof of~\ref{lem:KMboundsenergy}.

\textbf{Step 4. Smoothing.}
We extend $w$ by zero
to $\{x\le0\}$. Since $\operatorname{supp}\eta'\subset[\frac98, \frac{15}8]\subset\subset(1,2)$,
we can pick
a nonnegative mollifier
$\rho\in C^\infty_c(\R^2;[0,\infty))$, with integral one and supported in $B_{1/8}$,
and replace $w$ by $u:=\rho \ast w$.
Then $u\in C^\infty(\R^2)$, $u=0$ for $x\le1$,
so that~\ref{Th:KMzero}, \ref{lem:KMboundszero} and~\ref{enumcinfty}   hold. By convexity
$|\partial_2u|\le1$, so~\ref{Th:KMy1} holds.
The bounds in
\ref{Th:KMx23} and
\ref{lem:KMboundsdxu} also hold: They are automatic for $x<1$, and for $x\ge 1$ they survive the mollification with a worse constant,
since $x\ge1$ and $|x-\hat x|\le 1/8$ imply
$\hat x^{-1/3}\le 2 x^{-1/3}$
and
$\hat x^{2/3}\le 2 x^{2/3}$.

Finally we turn to the energy estimate~\ref{lem:KMboundsenergy}.
For $L\le 1$ the integral is still exactly $LH$, as above.
For $L>1$,
write $R:=(0,L)\times (a, a+H)$, and $\hat R:=
(0, L+1)\times (a-1, a+H+1)$.
By convexity,
\begin{equation*}
 \int_R |\partial_1 u|^2 +|D^2u|dx\,dy \le
 \int_{\hat R} |\partial_1 w|^2 +|D^2w|dx\,dy.
\end{equation*}
For the nonconvex term we need a finer estimate.
We first observe that for any $t,s\in[-1,1]$ one has
\begin{equation*}
 \min\{(1\pm s)^2\}
 \le \min\{(1\pm t)^2\}+4|t-s|
\end{equation*}
(it suffices to check that $t\mapsto (1-t)^2$ is $4$-Lipschitz in this interval), therefore it suffices to control $|\partial_2u-\partial_2w|$. In turn,
by the convolution estimate in $BV$ we have
\begin{equation}\label{eqpartialyuBVmoll}
 \int_R |\partial_2u-\partial_2w| dx\,dy
 \le C \int_{\hat R} |D^2w| dx\,dy .
\end{equation}
Combining the two bounds,
\begin{equation*}
E^1_\KM[u,R]\le C E^1_\KM[w,\hat R]\le C [(2L)^{1/3}(H+2)+L+1]
 \le C' (L^{1/3}H+L),
\end{equation*}
where in the final step we used $L\ge 1$.
This concludes the proof for
$\eps=1$, and therefore for all $\eps$.
\end{proof}

\subsubsection{Upper bound in a polyhedron} 
Let $P$ be a compact convex  polyhedron in $\R^n$ with nonempty interior, and let $F_1, \ldots, F_K$ denote its faces.
We use here $x$ for an element of $\R^n$.
We assume that
\begin{equation} \label{orientation_condition}
\text{ for any $k$, the vectors $\nu_k$ and $\bq$ are not parallel},
\end{equation}
{where $\nu_k$ is the inward unit normal to $F_k$.}
Note that the condition \eqref{orientation_condition}
is generically true, in the sense that, having fixed the polyhedron $P$  and the vector $\bq$, for almost any rotation $R\in O(n)$
the rotated polyhedron $RP$ has property~\eqref{orientation_condition}. Next, let us decompose $P = \cup_{k=1}^{K} P_k$, where $P_k$ is the polyhedron
\[
P_k:=\{x\in P: \dist(x, \d P)=\dist(x, F_k)\}.
\]
We shall control the full ``scalar'' energy, defined 
{for $\omega\subset\R^n$ open and bounded, $u\in W^{1,2}(\omega;\R^n)$ with $Du\in BV(\omega;\R^{n\times n})$} 
as
\begin{equation}\label{eqdefescalsect2}
\Eseps[u,\omega] := \int_{\omega} \bigl[\min|Du \pm \bp \otimes \bq|^2 + \ve |D^2 u|\bigr] \, dx
\end{equation}
Since  $|A_\sym|\le |A|$ for any matrix $A$,
recalling~\eqref{eq:energy}
one immediately obtains
$E_\ve[u,P]\le \Eseps[u,P]$.

\begin{lemma}\label{lem:UB2/3} Let 
$P\subset\R^n$ be a compact convex polyhedron with property~\eqref{orientation_condition}, $\bp$,  $\bq\in\R^n$. Then there is $C = C(\bp, \bq, P) > 0$ such that for any $\varepsilon>0$ there is $u$ as in~\eqref{eqfunctionspace} with
\begin{equation}\label{eq23polyhedr}
    \Eseps[u,P] \leq C \ve^{2/3}.
\end{equation}
Further
for any $L>0$,  $x_0\in\R^n$, there is $u$ as in~\eqref{eqfunctionspace}   with
\begin{equation}\label{eq23polyhedrL}
\Eseps[u, LP+x_0] \leq C(\bp, \bq, P) L^{n-2/3} \ve^{2/3}
\end{equation}
and $\mathrm{supp\,}u \subset\subset LP+x_0$.
    \end{lemma}
\begin{proof}
We start with~\eqref{eq23polyhedr}. If $\eps>1$, then $u=0$ will do, provided $C\ge |\bp|^2|\bq|^2|P|$. Assume $\eps\in(0,1]$, and
let  $u_\KM\in C^\infty((0,\infty)\times\R)$
be the function obtained from Lemma~\ref{lem:KMbounds}. We set
\begin{equation}\label{eqdefufromukmpolyh}
u(x):= \bp \, u_\KM(\dist(x,\d P), x \cdot \bq).
\end{equation}
In the following we shorten $d(x):= \dist(x, \d P)$.
Before computing the derivatives, we note that $d(x) $ is affine in each $P_k$, and thus $\nabla^2 d$ is a finite measure supported on $\cup_k\d P_k$, with
\begin{equation}\label{eqpropdP}
d\in W^{1,\infty}(P), \hskip5mm
|\nabla d|=1 \text{ a.e.}, \hskip5mm
Dd\in BV(P;\R^n), \hskip5mm
 |D^2d|\le 2\mathcal H^{n-1}\LL \cup_k\d P_k.
\end{equation}
To find the energy of $u$ in $P$, compute
\begin{equation}\label{first_der}
Du(x) = \d_1 u_\KM(d(x), x \cdot \bq)
\,
\bp \otimes \nabla d(x)
+ \d_2 u_\KM (d(x), x \cdot \bq)\,
\bp \otimes \bq
,
\end{equation}
as well as
\begin{multline}\label{second_der}
D^2 u (x) =  \d_1 u_\KM \, \bp\otimes \nabla^2 d(x) +\d_1^2 u_\KM\, \bp \otimes \nabla d(x) \otimes  \nabla d(x)
\\
+ \d_1 \d_2 u_\KM \, \bp\otimes (\nabla d(x) \otimes  \bq +\bq \otimes  \nabla d(x) )
+  \d_2^2 u_\KM\, \bp \otimes \bq \otimes \bq ,
\end{multline}
 where we write $u_\KM$ for $u_\KM(d(x), x \cdot \bq)$; $D^2u$ and $D^2d$ are understood as measures. We may now proceed with the estimate of the energy
treating separately the elastic and the surface contributions. The estimate for the elastic parts, in view of \eqref{first_der} and of $|\nabla  d(x)| = 1$ almost everywhere, reads as
\begin{multline}\label{elastic_bound}
\int_{P} \min \{\left|Du \pm (\bp \otimes \bq) \right|^2\}\, dx \\
\leq \int_{P} [2 |\bp|^2 |\d_1 u_\KM(d(x), x \cdot \bq)|^2 +
2 |\bp \otimes \bq|^2  \min \{(\d_2 u_\KM (d(x), x \cdot \bq)  \pm 1)^2\}] \, dx.
\end{multline}
Recalling~\eqref{eqpropdP} and that
$|\d_1 u_\KM|\le C$ almost everywhere
by~\ref{lem:KMboundsdxu}
and~\ref{lem:KMboundszero},
\begin{equation}\label{internal surface}
\ve \int_{P} \left| \bp \, \d_1 u_\KM \, \nabla^2 d(x) \right| dx \leq  C \ve |\bp| \sum_{k=1}^K \mathcal H^{n-1}(\d P_k)  \leq C(\bp, P) \ve.
\end{equation}
Thus, using \eqref{second_der} and \eqref{internal surface},
\begin{equation}\label{surface_bound}
\ve \int_{P}  |D^2 u| \, dx \leq \ve \int_{P} [|\bp| (1+|\bq|)^2 |D^2 u_\KM(d(x), x \cdot \bq)| ]\, dx + C(\bp, P) \ve.
\end{equation}
Combining \eqref{elastic_bound} and \eqref{surface_bound}, we obtain
\begin{equation*}
\Eseps[u,P] \leq C(\bp, \bq) \int_{P} \Psi(d(x), x \cdot \bq) \, dx + C(\bp,  P) \ve,
\end{equation*}
where we have set
\[
\Psi(t,y):=|\d_1 u_\KM(t,y)|^2 + \min\{(\d_2 u_\KM(t,y) \pm 1)^2 \} + \ve |D^2 u_\KM(t,y)|.
\]
We next use the coarea formula with respect to the map $\phi:P\to\R^2$, $\phi(x):=(d(x), x\cdot \bq)$.
Since  $\nabla d=\nu_k$ in $P_k$,
the Jacobian factor $|\nabla d\wedge \bq|$  takes finitely many values, which are nonzero
by \eqref{orientation_condition}.
We obtain
\begin{align*}
   & \int_{P} \Psi(d(x), x \cdot \bq) \, dx  = \int_{P} \frac{\Psi(d(x), x \cdot \bq)}{|\nabla d \wedge\bq|} \ |\nabla d \wedge \bq| \, dx \\
\nonumber =  &  \int_{0}^{\infty} \int_{\R} \int_{P\cap \{d(x) = t, x \cdot \bq =y\}} \frac{\Psi(t,y)}{|\nabla d \wedge\bq|}d \mathcal{H}^{n-2} \, dy \, dt \leq C(\bq, P) \int_{0}^{L} \int_a^{a+H} \Psi(t,y) \, dy dt
\end{align*}
with $L$ the diameter of $P$, $a:=-C|\bq|L$, $H:= 2C|\bq|L$.
We used that for any $t$ and $y$ one has
$\mathcal H^{n-2}(P\cap \{d(x) = t, x \cdot \bq =y\})\le C$, with $C$ depending only on $n$ and $P$.
In view of~\eqref{eq:KM-full-rectangle}, the proof of~\eqref{eq23polyhedr} is concluded.

Finally,~\eqref{eq23polyhedrL} follows by scaling. Indeed, set $\hat \eps:=\eps/L$, let $\hat u$ obey~\eqref{eq23polyhedr} in $P$
with $\hat \eps$, and define
\begin{equation}\label{eqscalepolyhedron}
  u(x_0+Ly):=L \hat u(y), \hskip1cm y\in P.
\end{equation}
By scaling and~\eqref{eq23polyhedr} $\Eseps[u,LP+x_0]
=L^n\Es_{\hat\eps}[\hat u, P]
\le C L^n (\eps/L)^{2/3}$, and the bound~\eqref{eq23polyhedrL} follows.
The condition on the support is immediate since
Lemma~\ref{lem:KMbounds}\ref{lem:KMboundszero} and~\eqref{eqdefufromukmpolyh} imply $u=0$ in a neighbourhood of $\partial P$.
\end{proof}
\subsubsection{Upper bound in a general domain} 
We are now ready to establish the upper bound in the first assertion of Theorem \ref{thm:main}.

\begin{proposition}\label{Le:2/3UB}
Let $\Omega\subset\R^n$ be as in Theorem \ref{thm:main}. Then there is $C>0$ such that for every $\ve>0$ we can find $u\in W_0^{1,2}(\Omega;\R^n)$ with $\nabla u\in BV(\Omega;\R^{n\times n})$ and
\[
\Eseps [u,\Omega]\le C \eps^{2/3}.
\]
\end{proposition}
The same statement, with the same proof, holds if one replaces the assumption that $\Omega$ is Lipschitz by the assumption that $\partial\Omega$ has bounded Minkowski content inside $\Omega$.
\begin{proof} If $\ve>1$, we can simply set $u \equiv 0$. The construction for $\ve \in (0,1]$ is based on a standard dyadic decomposition of $\Omega$ into disjoint cubes.
We will apply Lemma~\ref{lem:UB2/3} to each of these cubes, which requires~\eqref{orientation_condition} for their faces, so we cannot prescribe the orientation of the cubes independently of $\bq$. We fix $R\in SO(n)$ such that $\bq$ is not parallel to any of $Re_1,\dots,Re_n$,
and for $k\in\N$ we introduce the family of cubes
\[
\Pi_k:= \big\{ R(z+(0,2^{-k})^n):  z\in 2^{-k}\mathbb Z^n\big\},
\]
define
\[
\Omega_0: = \{ A \in \Pi_0: A \subset \Omega\},
\]
and then inductively
\[
\Omega_k: =\Bigl\{ A \in \Pi_k: A \subset \Omega \setminus \bigcup_{i=0}^{k-1}
\bigcup_{Q\in \Omega_{i}} Q \Bigr\}.
\]
By the dyadic structure, for any $k$, and any $A\in \Pi_k$ with $A\subset\Omega$, either 
$A\subset B$ for some $B\in \cup_{i<k} \Omega_i$ or $A\in \Omega_k$. Therefore
these cubes cover $\Omega$ up to a null set.
Similarly, if $Q\in\Omega_k$ 
and $k\ge 1$ then the unique $Q'\in \Pi_{k-1}$ with $Q\subset Q'$ is not contained in $\Omega$, so that $\dist(Q,\partial\Omega)\le C(n) 2^{-k}$. For $k=0$ one still has that the distance is bounded by the diameter of $\Omega$.

Since $\Omega$ is a Lipschitz set,
the Minkowski content of $\partial\Omega$ is bounded,
and we have
$|\{x\in \Omega: \dist(x,\partial\Omega)<r\}|\le C(\Omega)r$ for any $r>0$, which in particular implies
$\#\Omega_k \le C(\Omega) 2^{k(n-1)}$.

We define $k_0\in\N$ by
\[
\frac{1}{2^{k_0}} \leq \ve < \frac{1}{2^{k_0-1}}
\]
and let
\[
\Omega^*:=\bigcup_{i=0}^{k_0}
\bigcup_{Q\in \Omega_i}Q.
\]
At this point we define $u$. In each $Q\in \Omega_k$, $k\le k_0$,
$Q=R(z_Q+(0,2^{-k})^n)$,
we use the construction in Lemma~\ref{lem:UB2/3}
with $P=R[0,1]^n$, $L=2^{-k}$, $x_0=Rz_Q$.
In the rest we set $u=0$. Each cell function vanishes in a neighborhood of its cube boundary. Since only finitely many cells are used, their zero extensions satisfy~\eqref{eqfunctionspace}, with no additional gradient jumps at cube boundaries. The constant depends only on $P$, which in turn depends only on $R$ (and hence only on $\bq$).
By~\eqref{eq23polyhedrL},
\[
\Eseps [u, Q]
=
\Eseps[u,Rz_Q+2^{-k}R(0,1)^n]
\le C 2^{-k(n-2/3)} \eps^{2/3}.
\]
Summing over all those cubes,
\[
\Eseps [u, \Omega^*] \leq C \ve^{2/3} \sum_{k=0}^{k_0} \frac{2^{k(n-1)}}{2^{k(n-2/3)}} = C \ve^{2/3} \sum_{k=0}^{k_0} \frac{1}{(2^{1/3})^k} \leq C \ve^{2/3}.
\]
In turn, the boundary layer gives
\[\Eseps [u, \Omega\setminus \Omega^*]
=\Eseps [0, \Omega\setminus \Omega^*]
=|\bp|^2|\bq|^2|\Omega\setminus \Omega^*|\le C 2^{-k_0}
\le C\eps.
\]
Combining the two, and recalling that we are working with $\eps\le 1$, concludes the proof.
\end{proof}

\subsection{\texorpdfstring{Upper bound when $\bp$ and $\bq$ are parallel}{Upper bound when p and q are parallel}}

\subsubsection{Construction in a half-plane}

We first produce a construction in two dimensions, for the case that $\bp=\bq=e_2$. This is a linearization of the construction obtained in \cite{CC14,CC15}. Variants of this construction were used in \cite{RRT23,CKM24,MR26}.

In this section we use $(x,y)$ to denote points in the plane, as in Section~\ref{secKMhalfp}.
For $\ve>0$, $\omega\subset\R^2$ open and bounded, and
$v\in W^{1,2}(\omega;\R^2)$ with $Dv\in BV(\omega;\R^{2\times2})$, we write
\[
\ECCeps[v,\omega]
:=
\int_\omega
\left(
|\d_1v_1|^2
+
|\d_1v_2+\d_2v_1|^2
+
\min\{(\d_2v_2\pm1)^2\}
\right)\,dx\,dy
+
\ve |D^2v|(\omega).
\]

\begin{lemma}\label{eqconstrvcc}
There is $C>0$ such that for every $\ve>0$ there is a function
\[
v=v_\CC:(0,\infty)\times\R\to\R^2
\]
with the following properties:
\begin{enumerate}
\item\label{eqconstrvcczero}
$v_\CC(x,y)=0$ for $x\in(0,\ve)$, $y\in\R$;

\item\label{eqconstrvccCinfty}\label{eqconstrvccW12BV}
$v_\CC\in C^\infty((0,\infty)\times \R;\R^2)$;

\item\label{eqconstrvccenergy}
for every $L>0$, $H>0$, and $a\in\R$,
\begin{equation}\label{CC energy}
\ECCeps
[v_\CC,(0,L)\times(a,a+H)]
\le
C\left(
\ve^{4/5}L^{1/5}H+\ve L
\right);
\end{equation}

\item\label{eqconstrvccvcc1}
\begin{equation}\label{bound_v_1}
|v_{\CC,1}(x,y)|
\le
C\ve^{2/5}x^{3/5};
\end{equation}

\item\label{eqconstrvccvcc2}
\begin{equation}\label{bound_v_2}
|v_{\CC,2}(x,y)|
\le
C\ve^{1/5}x^{4/5};
\end{equation}

\item\label{eqconstrvccboundgrad}
\begin{equation}\label{bound_grad}
|Dv_\CC(x,y)|\le C.
\end{equation}
\end{enumerate}
\end{lemma}
\begin{proof}
{\bf Step 0. Rescaling.}
As in Theorem \ref{Th:KM} and Lemma \ref{lem:KMbounds}, we first reduce the statement to the
case $\ve=1$. Assume for the moment that there are $C>0$
and a function $\widehat v_\CC:(0,\infty)\times\R\to\R^2$
satisfying the conclusions of the lemma with $\ve=1$. For
$\ve>0$ define
\begin{equation}\label{CC_scaling}
v_\CC(x,y)
:=
\ve\,\widehat v_\CC\left(\frac{x}{\ve},\frac{y}{\ve}\right).
\end{equation}
Since $\widehat v_\CC$ has
properties~\ref{eqconstrvcczero}
and~\ref{eqconstrvccW12BV}, so does $v_\CC$.
Moreover, the pointwise estimates for $\ve=1$ give
\[
|v_{\CC,1}(x,y)|
\le
C\ve
\left(\frac{x}{\ve}\right)^{3/5}
=
C\ve^{2/5}x^{3/5},
\]
and
\[
|v_{\CC,2}(x,y)|
\le
C\ve
\left(\frac{x}{\ve}\right)^{4/5}
=
C\ve^{1/5}x^{4/5},
\]
so~\ref{eqconstrvccvcc1} and~\ref{eqconstrvccvcc2} also hold.
Further,
\[
Dv_\CC(x,y)
=
D\widehat v_\CC\left(\frac{x}{\ve},\frac{y}{\ve}\right),
\qquad \text{ and }
D^2v_\CC(x,y)
=
\frac1{\ve}
D^2\widehat v_\CC\left(\frac{x}{\ve},\frac{y}{\ve}\right).
\]
The first implies~\ref{eqconstrvccboundgrad}. Finally, changing variables $x=\ve\hat x$, $y=\ve\hat y$, we obtain
\begin{align}
&\ECCeps
[v_\CC,(0,L)\times(a,a+H)]
=
\ve^2
\ECCone
\left[
\widehat v_\CC,
\left(0,\frac{L}{\ve}\right)
\times
\left(\frac{a}{\ve},
      \frac{a+H}{\ve}\right)
\right].
\label{CC_energy_scaling}
\end{align}
Hence the estimate for $\ve=1$ implies
\begin{align*}
\ECCeps
[v_\CC,(0,L)\times(a,a+H)]
&\le
C\ve^2
\left[
\left(\frac{L}{\ve}\right)^{1/5}
\frac{H}{\ve}
+
\frac{L}{\ve}
\right]
=
C\left(
\ve^{4/5}L^{1/5}H+\ve L
\right),
\end{align*}
which is~\ref{eqconstrvccenergy}. Thus it suffices to prove the lemma for $\ve=1$.

\medskip
\noindent
{\bf Step 1. The period-doubling construction for $\ve=1$.}
We use the same period-doubling cell as in~\cite{CC15}.
Fix $l,h>0$, with $h\le 2l$, 
and choose a  function
$\gamma\in C^\infty(\R)$ such that
\[
\gamma(t)=0, \text{ for } t<\frac 18, \qquad
\gamma(t)=\frac18,\text{ for } t>\frac78,
\qquad \gamma'\ge0.
\]
This implies $|\gamma'|+|\gamma''|\le C$.
We set
\[
b_1(x):=\frac h8+h\gamma\bigl(\frac xl\bigr),
\qquad
b_2(x):=\frac{3h}{8}+h{\gamma\bigl(\frac xl\bigr)},
\]
\[
b_3(x):=\frac{5h}{8}-h{\gamma\bigl(\frac xl\bigr)},
\qquad
b_4(x):=\frac{7h}{8}-h{\gamma\bigl(\frac xl\bigr)}
\]
and then define $v^*_{l,h}:[0,l]\times\R\to\R^2$ by setting
\begin{equation}\label{eqdefvCCstar}
v^*_{l,h}(x,y):=
\begin{cases}
\begin{pmatrix}
0\\ y
\end{pmatrix},
&0\le y\le b_1(x),
\\[1.2em]
\begin{pmatrix}
\displaystyle
\frac{2h}{l}\gamma'(x/l)(-y+b_1(x))
\\[0.4em]
-y+2b_1(x)
\end{pmatrix},
&b_1(x)<y<b_2(x),
\\[1.2em]
\begin{pmatrix}
\displaystyle
-\frac{h^2}{2l}\gamma'(x/l)
\\[0.4em]
y-\frac h2
\end{pmatrix},
&b_2(x)\le y\le b_3(x),
\\[1.2em]
\begin{pmatrix}
\displaystyle
\frac{2h}{l}\gamma'(x/l)(y-b_4(x))
\\[0.4em]
\displaystyle
-y+\frac{3h}{4}-2h\gamma(x/l)
\end{pmatrix},
&b_3(x)<y<b_4(x),
\\[1.2em]
\begin{pmatrix}
0\\ y-h
\end{pmatrix},
&b_4(x)\le y\le h
\end{cases}
\end{equation}
and extending it periodically in the second variable
with period $h$, see Figure~\ref{fig:CC-cell}.
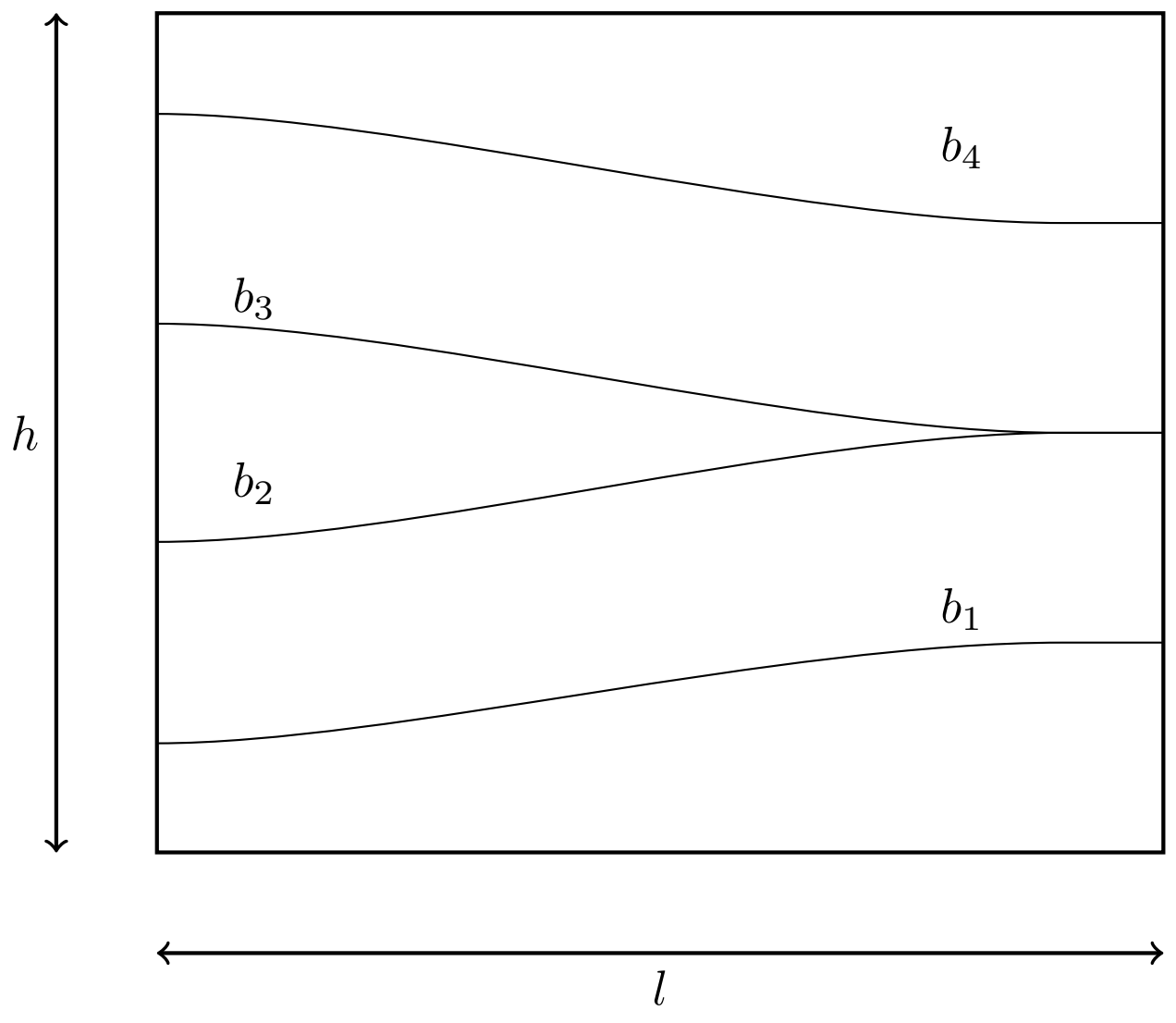
\begin{figure}[t]
\centering
\begin{tikzpicture}[scale=1.2]

\draw[thick] (0,0) rectangle (6,5);

\draw
(0,0.65)
.. controls (1.4,0.65) and (3.7,1.25) ..
(5.4,1.25)
-- (6,1.25);

\draw
(0,1.85)
.. controls (1.5,1.85) and (3.8,2.50) ..
(5.4,2.50)
-- (6,2.50);

\draw
(0,3.15)
.. controls (1.5,3.15) and (3.8,2.50) ..
(5.4,2.50)
-- (6,2.50);

\draw
(0,4.40)
.. controls (1.4,4.40) and (3.8,3.75) ..
(5.4,3.75)
-- (6,3.75);

\node at (4.8,4.20) {$b_4$};
\node at (0.58,3.3) {$b_3$};
\node at (0.58,2.2) {$b_2$};
\node at (4.8,1.45) {$b_1$};

\draw[<->,thick] (-0.6,0) -- (-0.6,5);
\node[left] at (-0.6,2.5) {$h$};

\draw[<->,thick] (0,-0.6) -- (6,-0.6);
\node[below] at (3,-0.6) {$l$};

\end{tikzpicture}
\caption{Branching building block in the $\bp = \bq $ case, see~\eqref{eqdefvCCstar}}
\label{fig:CC-cell}
\end{figure}
By construction $v^*_{l,h}$ is continuous and piecewise smooth; it obeys
\[
\d_1v^*_{l,h,2}+\d_2v^*_{l,h,1}=0,
\qquad
|\d_2v^*_{l,h,2}|=1
\quad\text{a.e.},
\qquad
{v^*_{l,h}(l,y)=v^*_{l,2h}(0,y)}
\]
and, since $h\le 2l$,
\begin{equation}\label{CCcellbounds-new}
|v^*_{l,h,1}|
\le C\frac{h^2}{l},
\qquad
|v^*_{l,h,2}|
\le Ch,
\qquad
|Dv^*_{l,h}|
\le
C. 
\end{equation}
These bounds imply
\begin{equation}\label{CCcellenergy-new}
\int_{(0,l)\times(0,h)}
\left(
|\d_1v^*_{l,h,1}|^2
+
|\d_1v^*_{l,h,2}+\d_2v^*_{l,h,1}|^2
+
\min\{(\d_2v^*_{l,h,2}\pm1)^2\}
\right)\,dx\,dy
\le
C\frac{h^5}{l^3}
\end{equation}
Differentiating the smooth pieces and adding the
gradient jumps along the four interface graphs gives,
after using $h\le2l$ to simplify the expression,
\begin{equation}\label{CCcellsurface-new}
 |D^2v^*_{l,h}|((0,l)\times(0,h))\le Cl.
\end{equation}
Indeed, the bulk contribution is bounded by
$C(h^2/l+h^3/l^2)$, while the jumps are bounded and the total
interface length is at most $Cl$. The periodic traces of the gradient
also agree at $y=0$ and $y=h$.

We now construct the self-similar pattern. As in the Kohn--M\"uller
construction, for $i\in\mathbb Z$ we set
\[
x_i:=\theta^i,\qquad
h_i:=\frac{1}{2^i},
\qquad
l_i:=x_i-x_{i+1}
=(1-\theta)x_i,
\]
where
\begin{equation}\label{choice_theta-new}
\theta:=2^{-5/4}\in\left(\frac14,\frac12\right)
\qquad\text{so that}\qquad
\theta^{4/5}=\frac12.
\end{equation}
On $(x_{i+1},x_i]\times\R$
we define
\[
\widetilde v_\CC(x,y)
:=
v^{\,*}_{l_i,h_i}
(x-x_{i+1},y).
\]
As in the original construction, the traces match at $x=x_i$,
because the trace at the right side of a cell of period $h_i$
coincides with the trace at the left side of a cell of period
$h_{i-1}=2h_i$. Hence $\widetilde v_\CC$ is well defined, continuous and piecewise smooth on
$(0,\infty)\times\R$.
Since $\gamma$ is constant near both endpoints, each cell is independent
of $x$ near its vertical sides, with first component zero. The matching
endpoint profiles therefore also have matching gradients almost
everywhere in $y$. There is no additional surface measure on $x=x_i$.

\medskip
\noindent
{\bf Step 2. Estimates for $\widetilde v_\CC$.}
By~\eqref{choice_theta-new},
\begin{equation}\label{eq2ixithetai}
2^{-i}=x_i^{4/5},
\qquad
(2\theta)^{-i}=x_i^{-1/5}.
\end{equation}
Using this and $\theta<1/2$,
\[
\frac{h_i}{l_i}
=
\frac{1}{1-\theta}(2\theta)^{-i}
\le
2x_i^{-1/5}.
\]
Therefore
\begin{equation}\label{CC_good_cells-new}
h_i\le 2l_i
\qquad\text{whenever}\qquad
x_i\ge 1, \qquad\text{which is the same as}
\qquad
i\le0.
\end{equation}
For $x\in[x_{i+1},x_i]$ we have $x_i\le x/\theta$. Thus,
using for each bound first~\eqref{CCcellbounds-new} and then~\eqref{eq2ixithetai},
\[
|\widetilde v_{\CC,2}(x,y)|
\le
Ch_i
=
C2^{-i}
\le
Cx^{4/5},
\]
and
\[
|\widetilde v_{\CC,1}(x,y)|
\le
C\frac{h_i^2}{l_i}
\le
C
\frac{4^{-i}}{\theta^i}
=C
\frac{x_i^{8/5}}{x_i}=C x_i^{3/5}
\le C x^ {3/5}.
\]
Finally, from~\eqref{CCcellbounds-new} and
\eqref{CC_good_cells-new},
\begin{equation}\label{CC_grad_aux-new}
|D\widetilde v_\CC(x,y)|
\le C
\qquad\text{for }x\ge 1.
\end{equation}

\medskip\newcommand\vCCt{w}
\noindent
{\bf Step 3. Modification close to the boundary.}
Let $\eta$ be the cutoff function introduced in~\eqref{eqchoiceeta} and define
\begin{equation}\label{def_vCC-new}
\vCCt(x,y)
:=
\widetilde v_\CC(x,y)
\eta(x).
\end{equation}
Then
\[
\vCCt(x,y)=0
\qquad\text{for }0<x<1.
\]
The preceding estimates immediately give
\[
|\vCCt_1(x,y)|
\le Cx^{3/5},
\qquad
|\vCCt_2(x,y)|
\le Cx^{4/5}.
\]
Furthermore,
\[
D\vCCt
=
\eta(x) D\widetilde v_\CC
+
\eta'(x)
\widetilde v_\CC\otimes e_1.
\]
Since $x\le 2$ in the transition region $[1,2]$, we have $|\widetilde v_\CC|\le C$ in this region, and therefore
\[
|D\vCCt|\le C.
\]
The construction also yields
\[
\vCCt\in
W^{1,2}_{\rm loc}((0,\infty)\times\R;\R^2),
\qquad
D\vCCt\in BV_{\rm loc}((0,\infty)\times\R;\R^{2\times 2}).
\]

\medskip
\noindent
{\bf Step 4. Energy estimate.}
We first estimate the energy of $\widetilde v_\CC$ away from the
cutoff region. Periodicity in $y$ and
\eqref{CCcellenergy-new}--\eqref{CCcellsurface-new} give, for $i\le0$,
\begin{equation}\label{CC_strip_energy-new}
\ECCone
[\widetilde v_\CC,
(x_{i+1},x_i]\times(a,a+H)]
\le
C\left(\frac{H}{h_i}+1\right)
\left(
\frac{h_i^5}{l_i^3}+l_i
\right)
\le
C\left(\frac{H}{h_i}+1\right)
l_i,
\end{equation}
where in the last step we used that the choice~\eqref{choice_theta-new} gives
$h_i^5(1-\theta)^4=l_i^4$.
Assume $L>1$ and let $I$ be the largest index with
$x_i\ge L$, which means that $L\in (x_{I+1}, x_I]$.
Summing the series,
\[
 \sum_{i=I}^0
C\left(\frac{H}{h_i}+1\right)
l_i
= C (x_I-x_1) +
CH \sum_{i=I}^0 {(2\theta)^i}\le C x_I + C H 2^{-I/4} \le C x_I + C H x_I^{1/5},
 \]
 where in the last step we used again~\eqref{choice_theta-new}.
Since $x_I\le L/\theta$,  we obtain for any $L\ge 1$
\[
\ECCone
[\widetilde v_\CC,
(1,L)\times(a,a+H)]
\le
C\left(L^{1/5}H+L\right).
\]
We turn to  the cutoff region. On
$(0,1)\times(a,a+H)$ we have $w=0$, and hence
\begin{equation}\label{eqenergCC01}
\ECCone
[w,(0,\hat L)\times(a,a+H)]
=\hat LH \le H\qquad\text{ for any $\hat L\in(0,1)$.}
\end{equation}
If $L\le1$, we are done. Otherwise,
on the transition strip $(1,2)\times(a,a+H)$,
from $|\widetilde v_\CC|+|D\widetilde v_\CC|
\le C$ we obtain
$|D\vCCt|\le C$, so that the elastic energy density is bounded,
and using
\[ |D^2\vCCt| \le
\eta |D^2\widetilde v_\CC|
+ C|D\widetilde v_\CC| + C|\widetilde v_\CC|
\]
we obtain
\[
\ECCone
[w,(0,2)\times(a,a+H)]
\le C H +
\ECCone
[\widetilde v_\CC,(1,2)\times(a,a+H)]
\le C (H+1).
\]
Since we are working in the case $L\ge 1$,
$H+1\le H L^{1/5}+L$.
Hence, for every $L,H>0$ and $a\in\R$,
\[
\ECCone
[\vCCt,(0,L)\times(a,a+H)]
\le
C\left(L^{1/5}H+L\right).
\]

\textbf{Step 5. Smoothing.}
We extend $w$ by zero to $x\le0$.
Use the same mollifier
$\rho\in C^\infty_c(\R^2;[0,\infty))$
as in Step 4 of
Lemma~\ref{lem:KMbounds}, and set $v_\CC:=\rho*w$.
Since $\rho$ is supported on $B_{1/8}$, and
$\eta=0$ on $(-\infty, 1+1/8)$, we have $v_\CC=0$ for $x\le1$.
Convolution preserves the uniform gradient bound and, because the
sampled positive first coordinates are comparable to $x\ge1$,
the componentwise bounds $|w_1|\le Cx^{3/5}$ and
$|w_2|\le Cx^{4/5}$.
As in Step 4 of
Lemma~\ref{lem:KMbounds}, the terms $(\partial_1 w_1)^2+(\partial_1w_2+\partial_2w_1)^2$ and
$|D^2w|$ are convex and therefore preserved by mollification (up to an irrelevant factor arising from the enlargement of the domain).
Finally, the same argument used in~\eqref{eqpartialyuBVmoll}
allows one to treat the term $\min|\partial_2w_2\pm 1|^2$.
The rescaling in Step 0 completes the proof.
\end{proof}

\subsubsection{Smooth distance to the boundary in a cube}

In order to proceed with the upper bound in a cube, we need to construct a smooth function that behaves like the distance to the boundary of a cube, such that its $k$-th order derivative is controlled by a constant times the distance to the boundary to the power $1-k$
(we shall only need a few values of $k$).
We focus on the reference cube
$(-1,1)^n$, and write
\begin{equation*}
\dist(x,\partial Q)=\min_i \min \{1\pm x_i\}, \hskip1cm
x\in  Q:=(-1,1)^n.
\end{equation*}
In order to make this function smooth we replace the minimum by the harmonic mean, and then estimate derivatives. A first property is that all
derivatives can be bounded by the appropriate negative power of the distance to the boundary, as specified in~\eqref{eqpropd1}.
However, for the change of variables it is crucial to have a lower bound on the first derivative, in a way that allows one to bound the Jacobian, see~\eqref{eqpropd2}. Since it vanishes at the origin, we shall then in Lemma~\ref{lemmacoarea} exploit polar coordinates to eliminate the singularity. For this reason we focus on obtaining a lower bound on the partial derivative in the radial direction.
\begin{lemma}\label{lemmaexistencedcube}
There are $C>0$ and a map $d\in C^\infty(Q;(0,\infty))$, both depending only on $n$, such that
\begin{equation}\label{eqpropd1}
    \frac1C \dist(x, \partial Q)\le d(x)\le \dist(x, \partial Q)
    , \quad
    |\nabla d|+ d |\nabla^2 d | + d^2|\nabla^3 d|\le C,
\end{equation}
and
\begin{equation}\label{eqpropd2}
    |x|^2\le -C x\cdot \nabla d (x).
\end{equation}
\end{lemma}

\begin{figure}\centering
 \includegraphics[width=5cm]{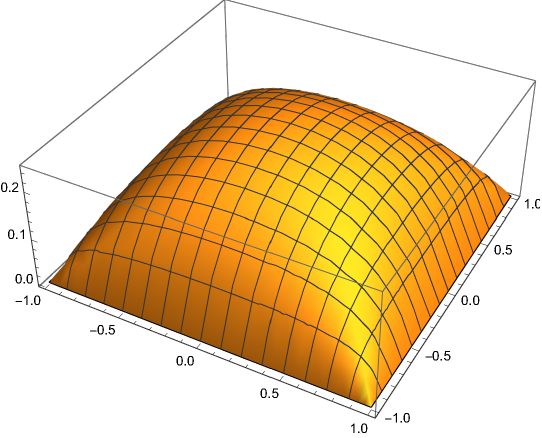}
 \caption{Plot of the function $d$ defined in~\eqref{eqdefgd}}
 \label{fig-plot-d}
\end{figure}

\begin{proof}
 We set, for $x\in Q$,
\begin{equation}\label{eqdefgd}
g(x):=\sum_{i=1}^n\Bigl( \frac1{1+x_i}+\frac1{1-x_i}\Bigr), \hskip6mm
 d(x):=\frac{1}{g(x)}=\left( \sum_{i=1}^n\Bigl( \frac1{1+x_i}+\frac1{1-x_i}\Bigr)\right)^{-1}
\end{equation}
(see Figure~\ref{fig-plot-d}).
Obviously $d\in C^\infty(Q)$,
and $d(x)\to0$ as $x\to\partial Q$.
We estimate the sum in $g(x)$ in terms of its largest element. An upper bound is obtained by multiplying by the number of elements,
\begin{equation}\label{eqboundgfromabove}
\begin{split}
g(x)
 &\le 2n \max_i \max \{
\frac1{1\pm x_i}\}
=\frac{2n}{\min_i \min\{1\pm x_i\}}
= \frac{2n}{\dist(x,\partial Q)},
\end{split}\end{equation}
so that $\frac1{2n}\dist(x,\partial Q)\le d(x)$.
Similarly,
\begin{equation*}
g(x) \ge  \max_i \max \{
\frac1{1\pm x_i}\}
= \frac{1}{\dist(x,\partial Q)},
\end{equation*}
so that
$d(x)\le\dist(x,\partial Q)$,
and the first bound in~\eqref{eqpropd1} holds for any $C\ge 2n$.

We next compute the partial derivatives,
\begin{equation}\label{eqggraddg}
\partial_id(x)=-\frac{\partial_i g(x)}{g^2(x)}=-d^2(x)\partial_i g(x), \hskip5mm
\partial_i g(x)=
\frac1{(1-x_i)^2}-\frac1{(1+x_i)^2}.
\end{equation}
Using $d(x)\le \dist(x,\partial Q)$ and $|\partial_i g|(x)\le \dist^{-2}(x,\partial Q)$ one sees that $Dd$ is bounded. Differentiating $Dd=-d^2Dg$,
\begin{equation*}
 D^2d=-2d Dg\otimes Dd-d^2D^2g
 = 2d^3 Dg\otimes Dg-d^2 D^2g,
\end{equation*}
and similarly
\begin{equation*}
|D^3d|\le 6d^4 |Dg|^3
+ 6 d^3 |D^2g|\, |Dg|+d^2|D^3g|.
\end{equation*}
From~\eqref{eqggraddg} one obtains $|D^kg|(x)\le C_k \dist^{-k-1}(x,\partial Q)\le C'_k d^{-k-1}(x)$ for any $k\ge 1$. Inserting this bound
into the last two equations gives the second bound in~\eqref{eqpropd1}.

To estimate the gradient from below we compute, again
from~\eqref{eqggraddg},
\begin{equation*}
 \partial_ig(x)=
 \frac{(1+x_i)^2-(1-x_i)^2}{(1+x_i)^2(1-x_i)^2}=
 \frac{4x_i}{(1-x_i^2)^2}
\end{equation*}
so that
\begin{equation*}
 x\cdot Dg(x)=4\sum_i \frac{x_i^2}{(1-x_i^2)^2},
 \hskip1cm
- x\cdot Dd(x)=4d^2(x)\sum_i \frac{x_i^2}{(1-x_i^2)^2}.
\end{equation*}
Pick an index $j$ with $|x_j|=\max_i|x_i|$, which is the same as $\dist(x,\partial Q)=\min\{1\pm x_j\}$.
Then $1-x_j^2=(1-x_j)(1+x_j)\le 2 \dist(x,\partial Q)$, and recalling
that  after~\eqref{eqboundgfromabove} we had shown
$d(x)\ge (2n)^{-1}\dist(x,\partial Q)$ we obtain
\begin{equation*}
 -x\cdot Dd(x)
 \ge  \frac{4d^2(x)}{(1-x_j^2)^2}x_j^2
 \ge\frac{ 4d^2(x) }{4\dist^2(x,\partial Q)}x_j^2
 \ge \frac{1}{(2n)^2} x_j^2 \ge \frac1{4 n^3} |x|^2,
\end{equation*}
which is~\eqref{eqpropd2}.
\end{proof}

%
%

We next show how this can be inserted into a change-of-variables formula.
\begin{lemma}\label{lemmacoarea}
Let $n\ge 2$, $M>0$, $C_1>0$, $d\in C^1(Q;[0,M])$ be such that
$-x\cdot Dd(x)\ge |x|^2/C_1$ pointwise. Then
there is $C=C(n,M, C_1)$ such that for any
Borel function $f:(0,M)\to[0,\infty)$ one has
 \begin{equation*}
  \int_Q f(d(x))dx
  \le C \int_0^M f(t)dt.
 \end{equation*}
\end{lemma}
\begin{proof}
We write, in polar coordinates,
\begin{equation*}
  \int_Q f(d(x))dx
  =\int_{S^{n-1}}
\int_{\{r>0: re\in Q\}} r^{n-1} f(d(re)) dr \,  d\calH^{n-1}(e).
\end{equation*}
For any fixed $e$, the map $r\mapsto d(re)$ is by  assumption $C^1$ and strictly decreasing, indeed, its derivative is
\begin{equation*}
\partial_r d(re)=e\cdot Dd(re)\le -r/C_1.
\end{equation*}
Therefore in the inner integral we can change variables from $r$ to $t:=d(re)$.
Since $n\ge 2$ and $r\le \sqrt n$,
\begin{equation*}
\begin{split}
\int_{\{r>0: re\in Q\}} r^{n-1} f(d(re)) dr
\le&
(\sqrt n)^{n-2} C_1
\int_{\{r>0: re\in Q\}} |\partial_r d(re)| f(d(re)) dr\\
\le &
(\sqrt n)^{n-2} C_1
\int_{(0,\max d)} f(t) dt.
\end{split}
\end{equation*}
This concludes the proof.
\end{proof}

\subsubsection{\texorpdfstring{Upper bound on subsets of $\R^n$}{Upper bound on subsets of R^n}}

We are finally in a position to establish the upper bound in the second assertion of Theorem~\ref{thm:main}. We recall that
\begin{equation*}
E_\eps [u,\Omega]=\int_{\Omega} \left(\min\{ \bigl| \frac{Du + Du^T}{2} \pm (\bp \otimes \bq) \bigr|^2 \}+\ve |D^2u|\right) \, dx .
 \end{equation*}

\begin{proposition}\label{Le:UBparallel}
Let $\Omega\subset\R^n$ be as in Theorem \ref{thm:main}, and suppose $\bp$ and $\bq$ are parallel. Then there is $C=C(\bp,\bq,\Omega)>0$ such that for all $\ve>0$ we can find $u\in W_0^{1,2}(\Omega;\R^n)$ with $\nabla u\in BV(\Omega;\R^{n\times n})$ and
\[
E_\eps [u,\Omega]\le C \eps^{4/5}.
\]
\end{proposition}
\begin{proof}
{\bf Step 1. We prove the assertion for $\Omega=Q:=(-1,1)^n$,
$\bp=\bq$.}

If $\ve>1$ the function $u \equiv 0$ has the required properties, so we can assume $\ve \in (0,1]$. Let
$v:[0,\infty)\times \R\to\R^2$ be the Chan--Conti construction from Lemma~\ref{eqconstrvcc},
and $d\in C^\infty(Q)$ be the smoothed distance function from Lemma~\ref{lemmaexistencedcube}.
For any translation $s\in(0,1)$ we define $u^s:Q\to\R^n$ by
\begin{equation}\label{eqdefuscc}
u^s(x) := \bp \, v_2(d(x), x \cdot \bp+s) + \nabla d(x) \, v_1(d(x), x \cdot \bp+s).
\end{equation}
The arguments of $v_1$ and $v_2$ are analogous to those in~\eqref{eqdefufromukmpolyh}. The necessity for the prefactor $Dd(x)$ in the second term, which is the one that requires control of third derivatives of $d$, is explained after~\eqref{eqstraindus}.
Implicitly evaluating $u^s$ and $d$ at $x$ and $v$ at $(d(x), x\cdot\bp+s)$,
\[
Du^s = (\bp \otimes \bp) \, \d_2 v_2 + (\bp \otimes \nabla d) \, \d_1 v_2 + (\nabla d \otimes \bp) \, \d_2 v_1 + (\nabla d \otimes \nabla d) \d_1 v_1 + (\nabla^2 d) \, v_1.
\]
Further,
\begin{align*}
| D^2 u^s| & \le | \bp|^2 |\nabla d | \, |\d_{21} v_2| 
+|\bp|^3 | \d_{22}v_2 | +|\bp| \, |\nabla d|^2  |\d_{11}v_2|+|\bp|^2| \nabla d|\,|\d_{12}v_2|  \\
    & + |\bp| |\nabla d|^2 |\d_{21} v_1 |+ 
    |\bp|^2|\nabla d|\,|\d_{22} v_1 | 
     + |\nabla d |^3 |\d_{11} v_1|+
     |\bp|\,|\nabla d |^2 |\d_{12} v_1 | \\
     & + |\bp|\,|\nabla^2 d| |\d_1 v_2|+ 2|\bp|\,|\nabla^2 d |\,| \d_2 v_1 |+ 3 |\nabla^2 d|\,|\nabla d| \,| \d_1 v_1 |\\
     & +
     | \nabla^3 d |\,|v_1|.
\end{align*}
The strain takes the form 
\begin{equation}\label{eqstraindus}
\frac{Du^s+(Du^s)^T}{2} = (\bp \otimes \bp) \, \d_2 v_2 + (\bp \otimes \nabla d+\nabla d \otimes \bp) \frac{ \d_1 v_2 + \d_2v_1}{2}
+  (\nabla d \otimes \nabla d) \d_1 v_1 + (\nabla^2 d) \, v_1 .
\end{equation}
This result motivates the choice in~\eqref{eqdefuscc}: the cancellation between $\partial_1 v_2$ and $\partial_2 v_1$ can again be exploited in higher dimension because the coefficients in~\eqref{eqdefuscc} match. The $Dd(x)$ which by the chain rule arises in the $\partial_1 v_2$ term  needs to be inserted manually in front of $v_1$ (the second term of~\eqref{eqdefuscc}). Since $|Dd|$ is bounded, there is
a constant $C$ that depends only on $n$ and $|\bp|$ with
\[
\left|\frac{Du^s+(Du^s)^T}2\pm \bp\otimes \bp\right|^2 \le C  (| \d_2 v_2\pm 1|^2 +|\d_1 v_2 + \d_2v_1|^2
+  |\d_1 v_1|^2) + 2|\nabla^2 d|^2 \, |v_1|^2 
\]
and
\[
|D^2 u^s| \leq C(|\nabla^2 v| + |\nabla^2 d| |\nabla v| + |\nabla^3 d| |v_1|).
\]
We separate the terms that compose the two-dimensional energy from the new error terms that arise from derivatives of $d$. Precisely, we write
\begin{equation}\label{eqEepsPsi12}
E_\eps[u^s,Q]
\leq C\int_{Q} \bigl(\Psi_1(d(x),x\cdot\bp+s)+\Psi_2(s,x) \bigr) \, dx
\end{equation}
where
$\Psi_1:(0,\infty)\times \R\to[0,\infty)$ is defined by
\begin{equation*}
\Psi_1:=
|\d_1 v_1|^2 +
|\d_1 v_2 + \d_2v_1|^2+
\min \{|\d_2 v_2\pm 1|^2\}+  \ve |D^2 v|,
\end{equation*}
and
$\Psi_2:(0,1)\times Q\to[0,\infty)$
by
\begin{equation*}
\Psi_2(s,x):=  |\nabla^2 d|^2\, |v_1|^2 + \ve |\nabla^2 d| |\nabla v| + \ve |\nabla^3 d| |v_1|
\end{equation*}
with $d$ evaluated at $x$ and $v$ evaluated at $(d(x),x\cdot \bp+s)$.

We start from the second term.
If $d(x)<\eps$, then $v=Dv=0$ and $\Psi_2(s,x)=0$.
If $d(x)\ge\eps$, \eqref{eqpropd1}, \eqref{bound_v_1}, and
\eqref{bound_grad} give,
writing again $d$ for $d(x)$,
\begin{equation*}
    |\Psi_2|(s,x)\le \frac{C}{d^2} \ve^{4/5} d^{6/5}+ C \ve \frac1{d} + C \ve \frac{1}{d^2} d^{3/5}\ve^{2/5}
    \le C \frac{\ve^{4/5}}{d^{4/5}} .
\end{equation*}
The bound is uniform in $s$. With Lemma~\ref{lemmacoarea}, recalling that $d(x)\le \dist(x,\partial Q)\le\sqrt n$,
\begin{equation}\label{eqboundpsi2}
\sup_{s\in (0,1)}
 \int_Q
 \Psi_2(s,x)dx
 \le C\eps^{4/5}
 \int_Q \frac{1}{d^{4/5}(x)} dx
 \le C \eps^{4/5} \int_0^{\sqrt n}
 \frac{1}{t^{4/5}}dt \le C \eps^{4/5}.
\end{equation}

We now turn to the first term.
Let $I:=(-|\bp|\sqrt n, 1+|\bp|\sqrt n)$.
This is defined so that
for any fixed $x\in Q$ one has
$x\cdot \bp+(0,1)\subseteq I$, and in particular, using $\Psi_1\ge0$,
\begin{equation*}
 \int_0^1\Psi_1(d(x),x\cdot\bp+s)\,ds
 \le
 \int_I\Psi_1(d(x),y)\,dy.
\end{equation*}
We integrate in $x$, swap the order of integration,
and use Lemma~\ref{lemmacoarea} as above,
\begin{equation*}\begin{split}
\int_0^1\int_{Q}\Psi_1(d(x),x\cdot\bp+s)\,dx\,ds
&\le\int_I \int_{Q}\Psi_1(d(x),y)\,dx\,dy
\le C \int_I \int_{0}^{\sqrt n} \Psi_1(t,y)\,dt\,dy.
\end{split}
\end{equation*}
By~\eqref{CC energy} with $L=\sqrt n$, $(a,a+H)=I$, the last integral is bounded by
$C\bigl(\eps^{4/5}|I|+\eps\bigr)
 \le C\eps^{4/5}$. Therefore we can choose $s\in(0,1)$ so that
\begin{equation*}
\int_{Q}\Psi_1(d(x),x\cdot\bp+s)\,dx
 \le C\eps^{4/5}.
\end{equation*}
Recalling~\eqref{eqEepsPsi12} and~\eqref{eqboundpsi2} the proof is concluded.

{\bf Step 2. Conclusion of the proof.}

We first prove the assertion for any axis-aligned cube, with the appropriate scaling. Precisely, we claim that  for any
 $x_0\in \R^n$, $L>0$,
there is $u$ with compact support in $x_0+(0,L)^n$ and
\begin{equation}\label{eqbdCCcubeq}
E_{\eps}[ u, x_0+(0,L)^n]
\le C L^n (\eps/L)^{4/5}.
\end{equation}
To do this, we scale by $L/2$ and translate.
Set $\hat \eps:=2\eps/L$.
By Step 1 there is $\hat u$ with compact support in $Q=(-1,1)^n$ and
$E_{\hat\eps} [\hat u,Q]\le C \hat\eps^{4/5}$.
Similarly to~\eqref{eqscalepolyhedron} we
let $x_*:=x_0+\sum_i \frac L2e_i$ and define
$  u(x_*+\frac L2y):=\frac L2 \hat u(y)$ for $y\in Q$;
by scaling we conclude $E_\eps[u,x_0+(0,L)^n]
=(\frac L2)^n E_{\hat\eps}[\hat u, Q]
\le C L^n (\eps/L)^{4/5}$, which is~\eqref{eqbdCCcubeq}.

The argument in the case of a general domain uses covering just like in the proof of \Cref{Le:2/3UB}, with the obvious adaptations in the exponents.

We also remark that the case when $\bp$ and $\bq$ are parallel reduces to the case when they are equal, since only the pair $\pm\bp\otimes \bq$ enters the energy.
\end{proof}

\section{Lower bounds}
\label{seclowerbound}

We next show how the lower bounds can be obtained by the localization method,  in the spirit of
\cite{Conti2006,CC15,CDMZ20}. The domain is covered by small cubes of side $\delta$; on each cube,
either the surface energy is large, or the elastic energy is large, or the displacement is large.
This can be understood as follows: if the surface energy is small, then the function is approximately affine; if also the elastic energy is small, then its gradient is close to one of the minima of the elastic energy. In turn, this implies that the function itself cannot be small. Details are given in
Lemma~\ref{lem:cube} below.
A global Poincar\'e-type estimate bounds the last alternative in
terms of the elastic energy, and optimizing in $\delta$ gives the scaling. The global estimate is
the only point where the two cases differ: if $\bp$ and $\bq$ are not parallel the elastic energy
controls the full displacement (Lemma~\ref{lem:poincare}), in the parallel case only with a weaker
exponent (Step~1 of the proof of Proposition~\ref{thm:LBparallel}).
Throughout this section, for $A\subset\R^n$ open we write
\[
 E_\ve[u,A]=E_{\rm el}[u,A]+\ve|D^2u|(A),\qquad
 E_{\rm el}[u,A]:=\int_A\min\{|(Du\pm\bp\otimes\bq)_\sym|^2\}\,dx .
\]
We recall that the three lower bounds presented below have been obtained by Fourier methods in~\cite{RRT23,RRTT24}, as discussed in the Introduction.

\subsection{Scaling and localization}

We start by discussing general aspects of the argument, which are common to the two cases.
Let $\omega\subset\R^n$ be nonempty and open.
We first rescale the domain and the vectors, to eliminate inessential notation, and then fix two global constants, $C_P$ and $c_*$.

\noindent\textbf{Scaling.}
Let $L>0$, $\bp=|\bp|\bp_0$, $\bq=|\bq|\bq_0$
be nonzero. For $u\in W^{1,2}_0(L\omega;\R^n)$ with $Du\in BV(L\omega;\R^{n\times n})$ we set
\begin{equation}
 u_0(y):=
 \frac{1}{L|\bp||\bq|}u(Ly),\qquad
 \ve_0:=\frac{\ve}{L|\bp||\bq|},
\end{equation}
so that $u_0\in W^{1,2}_0(\omega;\R^n)$ and $Du_0(y)=Du(Ly)/(|\bp||\bq|)$. A change of
variables gives
\begin{equation}\label{eqscalingLB2}
 E_\ve[u,L\omega]=L^n|\bp|^2|\bq|^2\,E^0_{\ve_0}[u_0,\omega],
\end{equation}
where $E^0$ is the energy with
$\bp_0$ and $\bq_0$ in place of $\bp$ and $\bq$. In particular, if
$E^0_{\ve_0}[u_0,\omega]\ge c\,\ve_0^\alpha$ for all $\ve_0\in(0,1]$, then
\begin{equation}\label{eqscalingLB3}
 E_\ve[u,L\omega]\ge c\,(|\bp||\bq|)^{2-\alpha}\ve^\alpha L^{n-\alpha}
 \qquad\text{for all }\ve\in(0,L|\bp||\bq|].
\end{equation}

\medskip
\noindent\textbf{Poincar\'e on cubes.}
We denote by $C_P$ a constant such that
\begin{equation}\label{eqBVPoinc}
 \int_{Q}|f-\bar f|\,dx\le C_P\,\delta\,|Df|(Q)
\end{equation}
for any cube $Q=x_0+(0,\delta)^n\subset\R^n$ of side $\delta$, any $f\in BV(Q;\R^m)$, and any $m\in\{1,n,n^2\}$, with $\bar f$
the average of $f$ on $Q$.  We also set
\[
 c_*:=\min_{|\eta|=1}\ \min_{b\in\R} \int_{(0,1)^n} |\eta \cdot x+b| \, dx > 0,
\]
so that, after translation and rescaling,
\begin{equation}\label{eqcstar}
 \int_{Q} |\eta \cdot x+b|\,dx\ge c_* \delta^{n+1}|\eta|
 \qquad\text{for all }\eta\in\R^n,\ b\in\R .
\end{equation}

The key ingredient
of the localization method is the fact that at a scale at which there is little surface energy, so that $u$ is approximately affine, there is a relation between the size of $u$ and the elastic energy. There are various ways of making this quantitative. We shall use the following lemma, which basically says that if both the surface energy and the $L^2$ norm of $u$ are small, then $Du$ is also small in $L^1$. It will become clear later that it is important to separate components; this is why we introduce the direction $\nu$.

\begin{lemma}[Localization on a cube]\label{lem:cube}
Let $Q=x_0+(0,\delta)^n\subset\R^n$ be a cube of side $\delta$, and let $u\in W^{1,1}(Q;\R^n)$ with
$Du\in BV(Q;\R^{n\times n})$. Then there is $\xi\in\R^{n\times n}$ such that
\begin{equation}\label{eqcubexi}
 \int_Q|Du-\xi|\,dx\le C_P\delta\,|D^2u|(Q)
\end{equation}
and, for every unit vector $\nu\in S^{n-1}$,
\begin{equation}\label{eqcubeB}
 c_*|\xi^T\nu|\,\delta^{n+1}\le
\delta^{n/2} \| \nu\cdot u\|_{L^2(Q)}+
C_P^2\delta^2|D^2u|(Q).
\end{equation}

\end{lemma}
\begin{proof}
Let $\xi$ be the average of $Du$ on $Q$; then \eqref{eqcubexi} is \eqref{eqBVPoinc} for $m=n^2$. Applying
\eqref{eqBVPoinc} with $m=n$ to $x\mapsto u(x)-\xi x$ we find $a\in\R^n$ with
\begin{equation}\label{eqpoincCp2}
 \int_Q|u-\xi x-a|\,dx\le C_P\delta\int_Q|Du-\xi|\,dx\le C_P^2\delta^2|D^2u|(Q).
\end{equation}
Since $\nu\cdot(\xi x+a)=(\xi^T\nu)\cdot x+\nu\cdot a$, \eqref{eqcstar} with $\eta=\xi^T\nu$ and the triangle inequality
give \begin{equation*}
 c_*|\xi^T\nu|\,\delta^{n+1}
 \le
 \int_Q|\nu\cdot u|\,dx
+ \int_Q|\nu\cdot (u-\xi x-a)|\,dx.
\end{equation*}
With H\"older's
inequality and~\eqref{eqpoincCp2} we then obtain~\eqref{eqcubeB}.
\end{proof}

\medskip
\noindent\textbf{Cubes.}
For $\delta>0$ let
\begin{equation}\label{eqdefAdeltacubes}
 A_\delta:=\{z\in \delta \mathbb Z^n: z+(0,\delta)^n\subset\omega\}, \qquad N_\delta:=\#A_\delta,
\end{equation}
and $Q_i:= x_i+(0,\delta)^n$, $i = 1,\ldots,N_\delta$.
As $\omega$ is open it contains a ball $B_r(x_0)$.
If $\delta$ is sufficiently small, on a scale set by $r$,
then every grid cube that meets $B_{r-\sqrt n\delta}(x_0)$ is contained in $\omega$, hence
\begin{equation}\label{estN}
 N_\delta\ge \frac{c_\omega}{\delta^n}
 \qquad\text{for all }\delta\in(0,\delta_\omega],\qquad
 \delta_\omega:=\frac{r}{2\sqrt n},\quad c_\omega:=|B_{r/2}|.
\end{equation}

\medskip
\noindent\textbf{Counting.}
Let $g$, $h\in L^1(\omega)$ be nonnegative (in practice,
$g$ will be the square of (a component
of) $u$ and $h$ the elastic energy density), and suppose that
we have chosen $C_A$, $C_B$, $C_C>0$ such that
each cube $Q_i$ satisfies at least
one of
\begin{equation}\label{eqABC}
 |D^2u|(Q_i)\ge C_A\delta^{n-1},\qquad
 \int_{Q_i}g\,dx\ge C_B\delta^{n+2},\qquad
  \int_{Q_i}h\,dx\ge C_C\delta^n .
\end{equation}
Then one of the three alternatives holds in at least $N_\delta/3$ cubes. Since the cubes are disjoint,
by \eqref{estN} for any $\delta\in(0,\delta_\omega]$ at least one of the following holds:
\begin{equation}\label{eqcounting}
 \ve|D^2u|(\omega)\ge\frac{c_\omega C_A}{3}\,\frac{\ve}{\delta},\qquad
 \int_\omega g\,dx\ge\frac{c_\omega C_B}{3}\,\delta^2,\qquad
 \int_\omega h\,dx\ge\frac{c_\omega C_C}{3}.
\end{equation}

\subsection{\texorpdfstring{Lower bound for $\bp$ and $\bq$ not parallel}{Lower bound for p and q not parallel}}

The global estimate is a Poincar\'e inequality which controls the $L^2$ norm of $u$ in terms of the elastic energy. This can be obtained using the boundary data and the one-dimensional Poincar\'e inequality in appropriate directions, chosen so that the term $\pm \bp\otimes \bq$ does not enter the estimates.
To explain the point, let us consider for example the case
 $\bp=e_1$, $\bq=e_1+e_2$. 
Then $(\bp\otimes \bq)_{\sym}=\begin{pmatrix}
    1 & 1/2\\ 1/2 & 0
\end{pmatrix}$ and  the elastic energy is
the integral of the density
\[
 \min_{\sigma\in\{-1,1\}}\left\{
 (\partial_1u_1+\sigma)^2+(\partial_2u_2)^2
 +\tfrac12(\partial_1u_2+\partial_2u_1+\sigma)^2\right\}.
\]
Therefore we can use Poincar\'e in the $x_2$ direction for $u_2$, but what to do for $u_1$?

\begin{lemma}\label{lem:poincare}
Let $\bp= |\bp|\bp_0$, $\bq= |\bq|\bq_0\in \R^n$ be nonzero and not parallel, $\omega\subset\R^n$ open and bounded, Lipschitz.
Then there is $C=C(\bp_0,\bq_0,{\rm diam\,}\omega)$ such that for any $u\in W^{1,2}_0(\omega;\R^n)$ one has
    \begin{equation*}
    \int_\omega |u|^2 dx \le C \int_\omega \min\{ |(Du\pm \bp\otimes \bq)_\sym|^2 \}dx.
\end{equation*}
\end{lemma}

\begin{proof}
Pick any vector $a\in S^{n-1}$. We write
\begin{equation*}
    a\cdot (Du\pm \bp\otimes \bq)_\sym a=
    a\cdot (Du\pm \bp\otimes \bq)a=
    a\cdot  Du \, a \pm (a\cdot \bp)(\bq\cdot a)
    .
\end{equation*}
If $a$ is orthogonal to either $\bp$ or $\bq$, then the second term vanishes.
 Using the one-dimensional Poincar\'e inequality  and the boundary data we have
\begin{equation*}\begin{split}
    \int_\omega (u\cdot a)^2  dx
    \le &C \int_\omega (\partial_a (u\cdot a))^2dx
    =C\int_\omega (a \cdot Du a)^2dx\\
    =&C\int _\omega(a \cdot (Du\pm \bp\otimes \bq)_\sym a)^2dx
    \le C\int_\omega \min\{|(Du\pm \bp \otimes \bq)_{\sym}|^2\}dx,
\end{split}\end{equation*}
with $C=({\rm diam\,}\omega)^2$.
Using that $\bp$ and $\bq$ are not parallel, we can find $n$ linearly independent vectors $a_1, \dots, a_n$
such that each of them is orthogonal to at least
one of $\bp$ and $\bq$. For example, one can first choose a basis $a_1, \dots, a_{n-1}$ of $\bp^\perp$ and then $a_n\in \bq^\perp\setminus \bp^\perp$. This exists since the two subspaces have dimension $n-1$ and are distinct by assumption.
We remark that this choice depends only on $\bp_0$ and $\bq_0$.
The vectors $\{a_i\}$ form a basis of $\R^n$, therefore there is $C>0$ such that $|\xi|^2\le C\sum_i|\xi\cdot a_i|^2$ for any $\xi\in \R^n$; summing over $i$ concludes the proof.
\end{proof}

As a warm-up, we start by establishing the lower bound for the energy functional
involving the full gradient $Du \pm \bp \otimes \bq$ rather than its symmetric part. We recall that in~\eqref{eqdefescalsect2} we had defined
\[
\Eseps[u,A] = \int_{A} \left(\min\{|Du \pm \bp \otimes \bq|^2 \} + \ve |D^2 u| \right) \, dx.
\]
While this result is auxiliary and not used in the following, it illustrates the idea of the lower bound in a less technical setting.

\begin{proposition}[Lower bound for the full gradient]\label{Th:LBnonSym}
Let $\omega\subset\R^n$ be nonempty, bounded, and Lipschitz, and let $\bp$, $\bq\in\R^n$ be nonparallel unit vectors.
Then there is $c = c(\omega, \bp, \bq)>0$ such that
\[
\Eseps[u,\omega] \geq c \min\{\ve^{2/3},1\}
\]
for every 
$u\in W^{1,2}_0(\omega;\R^n)$ with $Du\in BV(\omega;\R^{n\times n})$ and any $\eps>0$.
\end{proposition}
\begin{proof}
The problem is effectively scalar. The map
$v:= \bp \cdot u$ has $Dv=Du^T\bp$, so that $(Du\pm\bp\otimes\bq)^T\bp=Dv\pm\bq$, which implies
$|Dv\pm\bq|\le|Du\pm\bp\otimes\bq|$ for each sign, $|v|\le|u|$,  and $|D^2v|\le|D^2u|$.
 Lemma~\ref{lem:poincare} and $|A_\sym|\le|A|$ give
\begin{equation}\label{eqomegav2d2}
\int_\omega v^2\,dx\le\int_\omega|u|^2\,dx\le C\Eseps[u,\omega].
\end{equation}

We claim that there are positive constants $C_A$, $C_B$, $C_C$, which depend only on $n$, $C_P$ and $c_*$, such that for any $\delta\in(0,\delta_\omega]$ each cube satisfies at least one of
\begin{equation}\label{eqABC1}
  |D^2u|(Q_i)\ge C_A\delta^{n-1},\quad
  \int_{Q_i}v^2\,dx\ge C_B\delta^{n+2},\quad
  \int_{Q_i}\min\{|Du\pm\bp\otimes\bq|^2\}\,dx\ge C_C\delta^n ,
\end{equation}
which is \eqref{eqABC} with
$g:=v^2$ and $h:=\min\{|Du\pm\bp\otimes\bq|^2\}$.
The three alternatives have a clear meaning: either the surface energy is large, or the function itself is large, or the elastic energy is large. Once this is proven,
by~\eqref{eqcounting} and~\eqref{eqomegav2d2} we have
\[
 \Eseps[u,\omega] \geq c \min\left\{\frac{\ve}{\delta}, \delta^2, 1\right\}
 \qquad\text{for all }\delta\in(0,\delta_\omega].
\]
Choosing $\delta=\min\{\delta_\omega,\ve^{1/3}\}$ gives $\Eseps[u,\omega]\ge c\,\min\{\ve^{2/3},1\}$ and concludes the proof.

It remains to prove the claim. We fix the three constants so that
\begin{equation}\label{condition1}
C_B=(c_*/4)^2,\qquad
C_C^{1/2} + C_P C_A < \frac12
\qquad\text{and}\qquad
C_P^2 C_A \leq \frac{c_*}{4}.
\end{equation}
Assume that all three fail for a certain $\delta\in(0,\delta_\omega]$ and a certain $Q_i$, and let $\xi$ be as in Lemma~\ref{lem:cube}.
Then \eqref{eqcubeB} with $\nu=\bp$ gives
\begin{equation*}
 c_* |\xi^T\bp| \le C_B^{1/2} + C_P^2 C_A,
\end{equation*}
which by~\eqref{condition1} implies
$|\xi^T\bp|\le\frac12$.
Recalling $|\bq|=1$,
this implies $\min \{|\xi^T\bp \pm \bq|\} \geq \frac12$ and almost
everywhere in $Q_i$
\[
 \frac12 \leq \min\{|Dv \pm \bq|\} + |Dv - \xi^T\bp|\le\min\{|Du\pm\bp\otimes\bq|\}+|Du-\xi|.
\]
We integrate over $Q_i$, use H\"older's inequality, the failure of the third condition in~\eqref{eqABC1} and \eqref{eqcubexi}, and obtain
$\frac12\delta^n\le(C_C^{1/2}+C_PC_A)\delta^n$, which contradicts~\eqref{condition1}.
\end{proof}

This argument has been set up so that passing to the symmetrized gradient requires only a small change: the condition $|\xi^T\bp|\le \frac12$ is
replaced by $|\xi|\le 1/(2\sqrt2)$.

\begin{proposition}[Lower bound for $\bp$, $\bq$ not parallel]\label{Th:LBSym}
Let $\omega\subset\R^n$ be a nonempty, bounded Lipschitz set, $L>0$, $\Omega=L\omega$, and let
$\bp$, $\bq\in\R^n$ be nonzero and not parallel. Then there is $c = c(\omega, \bp/|\bp|, \bq/|\bq|)>0$ such
that
\[
 E_\ve[u,\Omega]\geq c\, (|\bp||\bq|)^{4/3}\ve^{2/3}L^{n-2/3}
\]
for every $\ve\in (0, L|\bp||\bq|]$ and every $u\in W^{1,2}_0(\Omega;\R^n)$ with
$Du\in BV(\Omega;\R^{n\times n})$.
\end{proposition}

\begin{proof}
By \eqref{eqscalingLB3} it suffices to show $E_\ve[u,\omega]\ge c\ve^{2/3}$
for $L=|\bp|=|\bq|=1$ and $\ve\in(0,1]$. Then
$|\bp\otimes \bq+\bq\otimes \bp|^2=2+2(\bp\cdot\bq)^2\ge 2$, that is,
$|(\bp\otimes \bq)_\sym|\ge 1/\sqrt2$. We set $\theta:=1/(2\sqrt{2n})$ and choose the three constants so that
\begin{equation}\label{condition1b}
 C_B=\Bigl(\frac{c_*\theta}{2}\Bigr)^2
,\qquad
 C_C^{1/2} + C_P C_A < \frac1{2\sqrt2},\qquad C_P^2 C_A \leq \frac{c_*\theta}{2}.
\end{equation}
We claim that for $\delta\in(0,\delta_\omega]$ each cube $Q_i$ satisfies one of the alternatives
\eqref{eqABC} with $g:=|u|^2$ and $h:=\min\{|(Du\pm\bp\otimes\bq)_\sym|^2\}$.

If this is true, then the conclusion follows as in the previous case: by \eqref{eqcounting} and Lemma~\ref{lem:poincare}, which gives
$\int_\omega|u|^2\,dx\le C E_\eps[u,\omega]$,
\[
 E_\ve[u,\omega] \geq c \min\left\{\frac{\ve}{\delta}, \delta^2, 1\right\}
 \qquad\text{for all }\delta\in(0,\delta_\omega].
\]
Choosing $\delta=\min\{\delta_\omega,\ve^{1/3}\}$ gives $E_\ve[u,\omega]\ge c\,\ve^{2/3}$ for
$\ve\in(0,1]$.

It remains to prove the claim.
Assume that all three fail for some $\delta$ and $Q_i$, and let $\xi$ be as in
Lemma~\ref{lem:cube}.
Then \eqref{eqcubeB} gives
\begin{equation*}
 c_* |\xi^T\nu| \le C_B^{1/2} + C_P^2 C_A \hskip5mm
 \text{ for any unit vector $\nu$},
\end{equation*}
which by~\eqref{condition1b} implies
$|\xi|\le\theta\sqrt n=1/(2\sqrt2)$.
In turn, this implies
$\min\{|(\xi \pm \bp \otimes \bq)_{\sym}|\}
\ge |(\bp\otimes \bq)_\sym|-|\xi|\geq \frac1{2\sqrt2}$, and almost everywhere in $Q_i$
\[
 \frac1{2\sqrt2} \leq \min\{|(Du \pm \bp \otimes \bq)_{\sym}|\} + |Du - \xi|.
\]
Integrating over $Q_i$, by H\"older's inequality, the failure of
the first and third alternatives in~\eqref{eqABC} and \eqref{eqcubexi},
\[
 \frac{1}{2\sqrt2} \delta^n \leq (C_C \delta^n)^{1/2}(\delta^n)^{1/2} + C_P C_A \delta^n,
\]
which contradicts~\eqref{condition1b}.
\end{proof}

\subsection{\texorpdfstring{Lower bound for $\bp$ and $\bq$ parallel}{Lower bound for p and q parallel}}

\begin{proposition}[Lower bound in the parallel case]\label{thm:LBparallel}
Let $\omega\subset\mathbb R^n$, $n\geq 2$, be bounded, Lipschitz and nonempty, let
$\Omega=L\omega$, $L>0$, and let $\bp$, $\bq\in\mathbb R^n$ be nonzero and parallel. Then there
exists $c=c(\omega)>0$ such that
\begin{equation}\label{eqlbpqparallmain}
 E_\varepsilon[u,\Omega] \geq c (|\bp||\bq|)^{6/5}\varepsilon^{4/5}L^{n-4/5}
\end{equation}
for every $\eps\in (0, L|\bp||\bq|]$ and every $u\in W^{1,2}_0(\Omega;\mathbb R^n)$ with
$Du\in BV(\Omega;\mathbb R^{n\times n})$.
\end{proposition}

\begin{proof}
Replacing the pair $(\bp,\bq)$ by $(\lambda \bp, \pm \lambda^{-1}\bq)$ for some $\lambda>0$ we can assume $\bp=\bq$;
this does not modify $|\bp||\bq|$ and $\{\pm \bp\otimes \bq\}$, and hence $E_\eps$.
By \eqref{eqscalingLB3} it suffices to show $E_\ve[u,\omega]\ge c\ve^{4/5}$ for $L=1$,
$|\bp|=|\bq|=1$, $\ve\in(0,1]$.
We further choose
$R\in SO(n)$ with $e_1=R\bp$,
and define $u_0:R\omega\to\R^n$ by
$u_0(y):= Ru(R^Ty)$, so that
\begin{equation*}
Du_0(y)\pm e_1\otimes e_1=RDu(R^Ty)R^T\pm (R\bp)\otimes (R\bp)=
R[ Du(R^Ty)\pm \bp\otimes \bp] R^T.
\end{equation*}
Therefore
we can further assume $\bp=e_1$,
if we can make sure that the constant is the same for all sets $R\omega$, $R\in SO(n)$. We write $u$ and $\omega$ for $u_0$ and $R\omega$ and work with
\begin{equation}\label{eqdefEele1e1}
 E_{\rm el}[u] =\int_\omega \min\left\{ |e(u)-e_1\otimes e_1|^2, |e(u)+e_1\otimes e_1|^2
 \right\}\,dx,
\end{equation}
where as usual $e(u)=(Du)_{\sym}$.

\medskip
\noindent
\textbf{Step 1: a coercivity estimate for the first component.}
The elastic energy does not control $u_1$ directly, but with a weaker exponent: we claim that
\begin{equation} \label{parallel-coercivity}
 \|u_1\|_{L^2(\omega)}^2 \leq C E_{\rm el}[u]^{1/2}
 \qquad\text{if } E_{\rm el}[u]\leq 1.
\end{equation}
To prove this, fix $j\in\{2,\ldots,n\}$ and
observe that  since
$u\in W^{1,2}_0(\omega;\mathbb R^n)$, the $2\times2$ minor formed by the first and the $j$-th rows
and columns is a null Lagrangian,
\begin{equation}\label{minor-zero}
 \int_\omega \left( \partial_1u_1\,\partial_ju_j - \partial_ju_1\,\partial_1u_j \right)\,dx =0
\end{equation}
(for $u\in C^\infty_c(\omega;\R^n)$ the integrand is
$\partial_1(u_1\partial_ju_j) -\partial_j(u_1\partial_1u_j)$; the general case follows by density).
Let $s_{1j}:=\frac12(\partial_1u_j-\partial_ju_1)$, so
that $\partial_ju_1=e_{1j}(u)-s_{1j}$ and $\partial_1u_j=e_{1j}(u)+s_{1j}$.
With $\partial_1u_1=e_{11}(u)$ and $\partial_ju_j=e_{jj}(u)$, \eqref{minor-zero} becomes
$\int_\omega ( e_{11}e_{jj} -e_{1j}^2 +s_{1j}^2 )\,dx=0$, and therefore
\begin{equation} \label{omega-est}
 \|s_{1j}\|_{L^2}^2 \leq \|e_{1j}\|_{L^2}^2 + \|e_{11}\|_{L^2}\|e_{jj}\|_{L^2}.
\end{equation}
From~\eqref{eqdefEele1e1} and $j\ne 1$ one obtains
$\|e_{1j}\|_{L^2}^2 + \|e_{jj}\|_{L^2}^2 \leq  E_{\rm el}[u]$. Moreover
$|e_{11}(u)| \leq 1+ \min \{|e_{11}(u)\pm 1|\}$ pointwise, hence $\|e_{11}\|_{L^2(\omega)}\leq C$
if $E_{\rm el}[u]\leq 1$. By \eqref{omega-est},
$\|s_{1j}\|_{L^2}^2 \leq C E_{\rm el}[u]^{1/2}$, and the one-dimensional Poincar\'e inequality in
the direction $e_j$ gives
\[
 \|u_1\|_{L^2(\omega)} \leq C\|\partial_ju_1\|_{L^2(\omega)}
 \leq C\|e_{1j}\|_{L^2}+C\|s_{1j}\|_{L^2} \leq C E_{\rm el}[u]^{1/4}.
\]
This proves \eqref{parallel-coercivity}, with a constant that only depends on the volume and diameter of $\omega$.

\medskip
\noindent
\textbf{Step 2: localization on cubes.}
We choose $C_A$, $C_B$, $C_C>0$ with
\[
 C_B=\Bigl(\frac{c_*}{4}\Bigr)^2
,\qquad
 C_PC_A+C_C^{1/2}<\frac12,\qquad C_P^2C_A\le\frac{c_*}{4},
\]
and claim that for $\delta\in(0,\delta_\omega]$ each cube $Q_i$ satisfies one of the alternatives
\eqref{eqABC} with $g:=|u_1|^2$ and $h:=\min\{|e(u)\pm e_1\otimes e_1|^2\}$. Assume that all three fail, and let $\xi$ be as in
Lemma~\ref{lem:cube}.
By~\eqref{eqcubeB} with $\nu=e_1$ we have
\begin{equation*}
 c_* |\xi_{11}|\le C_B^{1/2} + C_P^2 C_A,
\end{equation*}
so that $|\xi_{11}|\le 1/2$. In turn, $\min\{|\xi_{11}\pm1|\}\geq\frac12$, and almost
everywhere in $Q_i$
\[
 \frac12 \leq |\partial_1u_1-\xi_{11}| + \min\{|\partial_1u_1\pm1|\}
 \le |Du-\xi| + \min\{|e(u)\pm e_1\otimes e_1|\}.
\]
Integrating, by \eqref{eqcubexi}, H\"older's inequality and the failure of
the first and third alternatives in~\eqref{eqABC},
$\frac12\delta^n \leq C_PC_A\delta^n + C_C^{1/2}\delta^n$, a contradiction.

\medskip
\noindent
\textbf{Step 3: conclusion.}
By \eqref{eqcounting}, for every $\delta\in(0,\delta_\omega]$ either
$E_\ve[u,\omega]\ge c\min\{\ve/\delta,1\}$, or $\int_\omega|u_1|^2\,dx\geq c\delta^2$. In the
second case, either $E_{\rm el}[u]\ge1$, or \eqref{parallel-coercivity} gives
$E_\ve[u,\omega]\ge E_{\rm el}[u]\ge c\delta^4$. Hence
\[
 E_\varepsilon[u,\omega] \geq c\min\left\{ \frac{\varepsilon}{\delta}, \delta^4, 1 \right\}
 \qquad\text{for all }\delta\in(0,\delta_\omega],
\]
and $\delta:=\min\{\varepsilon^{1/5}, \delta_\omega\}$ gives $E_\varepsilon[u,\omega]\geq c\eps^{4/5}$
for $\ve\in(0,1]$; the constant only depends on the volume and diameter of $\omega$ and on the radius of the ball entering~\eqref{estN}.
\end{proof}

\begin{proof}[Proof of Theorem~\ref{thm:main}]
The upper bounds follow from \Cref{Le:2/3UB} and \Cref{Le:UBparallel}. \Cref{Th:LBSym}
and \Cref{thm:LBparallel} with $L=1$ give the lower bounds for $0<\eps\le\eps_*:=|\bp||\bq|$.
For $\alpha=2/3$ or $4/5$, respectively, monotonicity of $E_\eps[u,\Omega]$ in $\eps$ gives, when
$\eps_*\le\eps\le1$,
\[
 E_\eps[u,\Omega]\ge E_{\eps_*}[u,\Omega] \ge c\eps_*^\alpha\ge c\eps_*^\alpha\eps^\alpha.
\]
Adjusting the constants proves the bounds on the entire interval $(0,1]$.
\end{proof}

\section*{Use of artificial intelligence}
The results and a complete first draft of this paper were obtained without generative AI.
The authors subsequently used Claude (Anthropic, Claude Opus~5.5) to check proofs, to improve
the presentation, and to assist with the literature review. All changes were made, checked
and approved by the authors, who take full responsibility for the content.
 
\section*{Acknowledgements}
We are grateful to the late Robert V. Kohn for his support and guidance.
This work was partially supported by the Deutsche Forschungsgemeinschaft (DFG) through project CRC1720 -- 539309657. The research of Oleksandr Misiats was supported by the Simons Foundation through Collaboration Grant for Mathematicians No. 854856 and National Science Foundation Grant DMS-2408507.

\section*{Statements and Declarations}
Funding information is provided in the Acknowledgements section. The authors declare that they have no competing interests. No datasets were generated or analysed during this study.

\begingroup

\endgroup

\end{document}